\documentclass[11pt,twoside]{article}
\usepackage{fullpage}
\usepackage{amssymb}
\usepackage{amsmath}
\usepackage{epsfig}

\def \ovA{${\bar{\mathbf A}}$}

\def \cF{\mathcal F}  
  
\def \cC{\mathcal C}  
\def \bR{\mathbb R}  
\def \bZ{\mathbb Z}

\def \tz0{\tau^Z_0}  
  
\def \tzp0{\tau^{Z'}_0}  
  
\newcommand{\bbN}{{\mathbb N}}
\newcommand{\bbR}{{\mathbb R}}
\newcommand{\bbE}{{\mathbb E}}

\newcommand{\bbZ}{{\mathbb Z}}

\newcommand{\mysection}{\setcounter{equation}{0} \section}  

\newtheorem{thm}{Theorem}[section]
\newtheorem{prop}[thm]{Proposition}
\newtheorem{cor}[thm]{Corollary}
\newtheorem{lem}[thm]{Lemma}

\begin{document}

\parindent=0pt

\thispagestyle{empty}

\title{Distributed Algorithms in an Ergodic  Markovian Environment}

\author{Francis Comets${}^{(a),}$\footnote{comets@math.jussieu.fr; http://www.proba.jussieu.fr/pageperso/comets/comets.html} \hspace*{.05cm},
Fran\c{c}ois Delarue${}^{(a),}$\footnote{delarue@math.jussieu.fr; http://www.math.jussieu.fr/$\sim$delarue} \hspace*{.05cm} and 
Ren\'e Schott${}^{(b),}$\footnote{schott@loria.fr; http://www.loria.fr/$\sim$schott}
 \\ \\ (a) Laboratoire de Probabilit\'es et Mod\`eles
Al\'eatoires,\\ Universit\'e Paris 7, UFR de Math\'ematiques, Case 7012, \\ 2, Place Jussieu, 75251 Paris Cedex 05 - FRANCE.
\\
\\
(b) IECN and LORIA, Universit\'e Henri Poincar\'e-Nancy 1, \\ 54506 Vandoeuvre-l\`es-Nancy - FRANCE.} 
\maketitle
\begin{abstract}
We provide a probabilistic analysis of the $d$-dimensional banker algorithm when transition probabilities may depend on time and space.
The transition probabilities evolve, as time goes by, along the trajectory of an ergodic Markovian environment, whereas the 
spatial parameter just acts on long runs. Our model complements the one considered by Guillotin-Plantard and Schott \cite{guillotin schott 2002} 
where transitions are governed by a dynamical system, and appears as a new (small) step towards more general time and 
space dependent protocols. 
\vspace{5pt}
\\
Our analysis relies on well-known results in stochastic homogenization theory and investigates the asymptotic behaviour
of the rescaled algorithm as the total amount of resource available for allocation tends to the infinity. In the two dimensional setting,
we manage to exhibit three different possible regimes for the deadlock time of the limit system. To the best of our knowledge, 
the way we distinguish these regimes is completely new.
\vspace{5pt}
\\
\textit{
\vspace{5pt}
}
\end{abstract}
\mysection{Introduction}
\label{Section 1}
Many real-world 
phenomena involve time and space dependency.
Think about option pricing: the behavior of traders is not the same when stock markets are opening and a few minutes before closure. Multi-agents problems, dam management problems are typical examples where (often random) decisions have to be made under time and space constraints.
\vspace{5pt}
\\
In computer science 
such
problems are usually called resource sharing problems.
Their statement is as follows (see Maier \cite{maier 1991}, Maier and Schott \cite{maier schott 1993}):
\vspace{5pt}
\\
Consider the interaction of $q$~independent processes
$P_1,\dots,P_q$, each with its own memory needs. The processes are allowed to
allocate and deallocate $r$~different, non-substitutable resources (types
of memory): $R_1,\dots,R_r$.  Resource limitations and 
resource exhaustions are defined as~follows. At~any time~$s$, process~$P_i$ is assumed
to have allocated some quantity $y_i^j(s)$ of resource~$R_j$  (both
time and resource usage are taken to be discrete, so~that $s\in\bbN$ and
$y_i^j(s)\in\bbN$).  Process~$P_i$ is assumed to have some maximum
need~$m_{ij}$ of resource~$R_j$, so~that
\begin{equation}
\label{eq:needconstraint}
0\le y_i^j(s) \le m_{ij}
\end{equation}
for all~$s$. The numbers 
$m_{ij}$~may be infinite; if~finite, it~is a hard limit that
the process~$P_i$ never attempts to exceed.  The resources~$R_j$ are
limited, so~that
\begin{equation}
\label{eq:resconstraint}
\sum_{i=1}^q y_i^j(s)< \Lambda_j
\end{equation}
for $\Lambda_j-1$ the total amount of resource~$R_j$ available for allocation.
Resource exhaustion occurs when some process~$P_i$ issues an unfulfillable
request for a quantity of some resource~$R_j$.  Here `unfulfillable' means
that fulfilling the request would violate one of the
inequalities~(\ref{eq:resconstraint}).
\vspace{5pt}
\\
The state space~$Q$ of the memory allocation system is the subset 
of~$\bbN^{qr}$ determined by~(\ref{eq:needconstraint})
and~(\ref{eq:resconstraint}).  This polyhedral state space is familiar:
it~is used in the banker algorithm for deadlock avoidance. Most
treatments of deadlocks (see Habermann \cite{haberman 1978})
assume that processes request and release resources in a mechanical way:
a~process $P_i$ requests increasing amounts of each resource~$R_j$ until
the corresponding goal~$m_{ij}$ is reached, then releases resource units
until $y_i^j=0$, and repeats  (the $r$~different goals of the process need
not be reached simultaneously, of~course).  This is a powerful assumption:
it~facilitates a classification of system states into `safe' and `unsafe'
states, the latter being those which can lead to deadlock.  However it is
an idealization.
Assume that regardless
of the system state, each process~$P_i$ with $0<y_i^j<m_{ij}$ can issue
either an allocation or deallocation request for resource~$R_j$.  The
probabilities of the different sorts of request may depend on the current
state vector~$(y_i^j)$.  In~other words the state of the storage
allocation system is taken as a function of time to be a finite-state Markov chain;
this is an alternative approach which goes  back at least as far as
Ellis \cite{ellis 1977}.
\vspace{5pt}
\\
The goal is to estimate the amount of time~$\tau$ until memory
exhaustion occurs, if~initially the $r$~types of resources are completely
unallocated: $y_i^j=0$ for all~$i,j$.  
The consequences of expanding the resource
limits~$\Lambda_j-1$ and the per-process maximum needs~$m_{ij}$ (if~finite)
on~the expected time to exhaustion are particularly interesting for practical applications.
There has been little work on the exhaustion of {\em shared\/} memory, or on `multidimensional' exhaustion, where one of a number of 
inequivalent resources becomes exhausted.
D. E. Knuth \cite{knuth 1973}, A. Yao \cite{yao 1981}, P. Flajolet \cite{flajolet 1986},  G. Louchard and R. Schott \cite{louchard schott 1991}, have provided
combinatorial or probabilistic analysis of some resource sharing problems under the assumption that transition probabilities are constant. Maier provided a large deviation analysis of colliding stacks for the more difficult case in which the transition probabilities are nontrivially state-dependent. More recently
N. Guillotin-Plantard and R. Schott \cite{guillotin schott 2002} analysed a model of exhaustion of shared resources where allocation and deallocation requests are modelled by time dependent dynamic random variables. 
\vspace{5pt}
\\
In this paper, we analyse such problems when the 
probability transitions are time and space dependent. 
We incorporate in the transitions of our model, the influence of 
an environment which evolves randomly albeit in a stationary manner
on the long run.
Our analysis relies 
on well-known results in stochastic homogenization. 
We focus on the regime when memory exhaustion is caused by normal
fluctuations (as in many of the above references), but
not by large deviations (as in Maier \cite{maier 1991}).
This paper can be viewed as a (small) step towards the analysis of protocols where decisions are time and space dependent random variables. 
In addition to incorporating space dependence,  we improve on
the results of Guillotin-Plantard and Schott \cite{guillotin schott 2002} in covering larger 
perturbations of time-homogeneous random walks.
\vspace{5pt}
\\
The organization is as follows: Section  \ref{banker} presents an example of resource sharing algorithm usually called banker algorithm. 
Precise statements of the main results are given in Section \ref{Section_Main_Results}. Proofs are postponed to the three 
following parts of the paper:
we investigate the asymptotic behaviour of the
algorithm in Sections \ref{ProofTHM1} and \ref{Proof_THM2} and we discuss in Section \ref{behaviour} several different regimes for the limit
process. As a conclusion, we indicate 
further conceivable extensions in Section \ref{Conclusion}.
\mysection{Example: the Banker Algorithm}
\label{banker}
\subsection{Practical Description of the Model in Dimension Two}
\label{banker.I}
We consider a simple distributed algorithm which involves only two customers $C_1$ and $C_2$ sharing a fixed quantity of a given resource M (money). 
There are fixed upper bounds $m_1$ and $m_2$ on how much of the resource each of the customers is allowed to use at any time. 
The banker decides to give to the customer $C_i$, $i=1,2$, the required resource units only if the remaining units are sufficient in order 
to fulfill the requirements of $C_j$, $j=1,2;j\neq i$.
This situation is modelled 
(see  Figure 1) 
by a random walk in a rectangle \index{rectangle} with a broken corner 
\index{broken corner}, i.e.,\\ 
\[
\{(x_1,x_2) \in \bbZ^d, \ 0\leq x_1\leq m_1, \ 0\leq x_2\leq m_2, \ x_1+x_2\leq \Lambda  \}
\] 
where the last constraint generates the broken corner. 
The random walk is reflected on the sides parallel to the axes and is absorbed on the sloping side.

\begin{figure} [htb]
\begin{center}
\includegraphics[
width=0.3\textwidth,angle=0]
{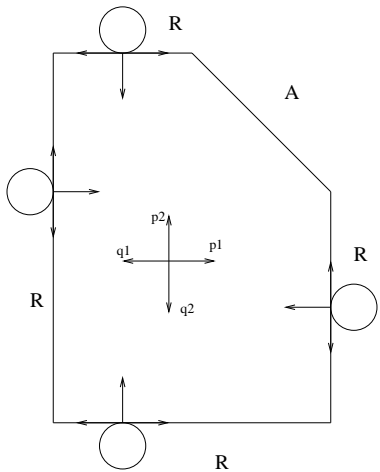}
\caption{Banker algorithm}
\label{fig:bankeralgorithm}
\end{center}
\end{figure}
Again the hitting place and the hitting time of the absorbing boundary are the parameters of interest.
This algorithm has been analysed by Louchard et al.\cite{louchard 1995, louchard schott 1991, louchard tolley zimmerman 1994}
in the model where the transitions are constant
and by Guillotin-Plantard and Schott \cite{guillotin schott 2002} when the transitions are some time dependent dynamic random variables.
\vspace{5pt}
\\ 
Our goal here is to analyse such distributed algorithms in an ergodic markovian environment where the transition probabilities are time and space dependent,
in a sense that we specify below. The environment takes care of purely unpredictable behaviors
in addition to periodic and quasi-periodic ones. 
\vspace{5pt}
\\
For the sake of simplicity, we assume from now on that $m_1=m_2=m$ and that $\Lambda$ writes $\Lambda = \lambda m$, with $\lambda \in [1,2)$.
The size of the problem is then given by the parameter $m$, which therefore plays a crucial role in
the time scale. For large $m$'s, we will prove approximations 
for the deadlock time by taking advantage of averaging properties of
the environment.    
\subsection{Probabilistic Modelization}
\label{banker.II}
We generalize the situation described in Subsection \ref{banker.I} and assume that there are $d$ customers instead of two. We then modelize the resources
of these customers by a Markov chain with values in the broken corner 
$\{(x_1,\dots,x_d) \in \bbZ^d:0\leq x_i\leq m, 1 \leq i \leq d, \ x_1+\dots+x_d \leq \lambda m  \}$, $m \in \bbN^*$ and $\lambda \in [1,d)$.
\vspace{5pt}
\\
{\bf \emph{Walk without Boundary Conditions.}}
\vspace{5pt}
\\
Define first the transitions of the chain without taking care of the boundary conditions. Recall to this end 
that the transition matrix of a time-space homogeneous random walk to the 
nearest neighbours in ${\mathbb Z}^d$, $d \geq 1$,
 reduces to a probability $p(\cdot)$ on the $2d$ directions of the 
discrete grid,
$${\mathcal V} \equiv \{e_1,-e_1,\dots,e_d,-e_d\}\;,$$
where $(e_i)_{1 \leq i \leq d}$ denotes the canonical basis of ${\mathbb R}^d$.
For given $i \in \{1,\dots,d\}$ and $u \in {\mathcal V}$, $p(u)$ simply denotes the 
probability of going from the current position $x$ to $x+u$.
\vspace{5pt}
\\
Assume for a moment that the transition matrices are
 space homogeneous but depend on time through
some 
environment (evolving with time)
with values in a finite space $E$, $N\equiv |E|$. 
We then need to consider, not a single transition probability, but
 a family
$p(1,\cdot)$, \ldots, $p(N,\cdot)$  of $N$
probabilities on ${\mathcal V}$. 
If the environment at time $n$ is in the state $i \in E$, then the transition
of the Markov chain at that time is governed by the probability $p(i,\cdot)$.
\vspace{5pt}
\\
When the environment is given by  a stochastic process $(\xi_n)_{n \geq 0}$
on $E$, the jumps  $(J_n)_{n \geq 0}$ of the random walk
are such that, for every $u \in {\mathcal V}$:
\begin{equation*}
{\mathbb P}\bigl\{J_{n+1}=u|{\mathcal F}_n^{\xi,J} \bigr\} =p(\xi_n,u),
\end{equation*}
with ${\mathcal F}^{\xi,J}_n \equiv \sigma \{\xi_0,\dots,\xi_n,J_0,\dots,J_n\}$.
\vspace{5pt}
\\
From now on, we assume that the environment $(\xi_n)_{n \geq 1}$ is a  
time-homogeneous Markov
chain on $E$, and we denote by  $P$ its transition matrix,
$P(k,\ell) \equiv {\mathbb P}(\xi_{n+1}=\ell|\xi_{n}=k)$ for $ k,\ell \in E$.
Then, the couple $(\xi_n,J_n)_{n \geq 1}$ is itself a time-homogeneous
Markov chain on the product space
$E \times {\mathbb Z}^d$ governed by the following transition:
\begin{equation*}
\forall k \in E, \ \forall u \in {\mathcal V}, \ {\mathbb P} \bigl\{\xi_{n+1}= k ,J_{n+1}=u |{\mathcal F}^{\xi,J}_n\bigr\}
= P(\xi_{n},k) p(\xi_{n},u).
\end{equation*}
Define the position of the walker in ${\mathbb Z}^d$:
\begin{equation} \label{183}
S_0 \equiv 0, \
\quad \forall n \geq 0, \ S_{n+1} \equiv S_n + J_{n+1}.
\end{equation}
In view of the applications mentioned above, the model is not fine enough. For this reason, 
we also assume that the steps $(J_n)_{n \geq 1}$ depend on the walker position in the following way:
for all  $k \in E $ and $u \in {\mathcal V}$,
\begin{equation}
\label{star}
{\mathbb P} \bigl\{\xi_{n+1}=k,J_{n+1}=u |{\mathcal F}^{\xi,J}_n\bigr\}
= P(\xi_{n},k) p(\xi_{n},S_n/m,u),
\end{equation}
where $m$ denotes a large integer that refers to the size of the box in Subsection \ref{banker.I} and, for each $k \in E$ and
$y \in \bbR^d$, $p(k,y,\cdot)$ a probability on 
${\mathcal V}$. 
Note that the random walk $(S_n=S_n^{(m)})_{n \geq 0}$ depends on the parameter 
$m$. Nevertheless, for simplicity we will often forget the dependence on
$m$ in 
our notations. In other words, $(\xi_n,S_n)_{n \geq 0}$ defines a Markov chain with rates:
\begin{equation*}
\forall k \in E, \ \forall u \in {\mathcal V}, \ {\mathbb P} \bigl\{\xi_{n+1}=k,S_{n+1}=u+S_n |{\mathcal F}^{\xi,S}_n\bigr\}
= P(\xi_{n},k) p(\xi_{n},S_n/m,u),
\end{equation*}
where ${\mathcal F}^{\xi,S}_n = \sigma \{\xi_0,\dots,\xi_n,S_0,\dots,S_n\}$.
\vspace{5pt}
\\
{\bf \emph{Walk with Reflection Conditions.}}
\vspace{5pt}
\\
Our original problem with reflection on the hyperplanes
$y_i=0, \ i \in \{1,\ldots,d\}$ and $y_i=m, \ i \in \{1,\ldots, d\}$ follows from
a slight correction of the former one. The underlying reflected walk $(R_n)_{n\geq 0}$
(also denoted by $(R_n^{(m)})_{n \geq 0}$ to specify the dependence on $m$)
satisfies, with 
${\mathcal F}^{\xi,R}_n \equiv \sigma \{\xi_0,\dots,\xi_n,R_0,\dots,R_n\}$:
\begin{equation*}
\forall k \in E, \ \forall u \in {\mathcal V}, \ {\mathbb P} \bigl\{\xi_{n+1}=k,R_{n+1}=u+R_n |{\mathcal F}^{\xi,R}_n \bigr\}
= P(\xi_{n},k) q(\xi_{n},R_n/m,u),
\end{equation*}
where $q$ denotes the kernel:
\begin{equation*}
\forall k \in E, \ \forall y \in [0,1]^d, \ \forall \ell \in \{1,\dots,d\}, \ q(k,y,\pm e_{\ell}) = p(k,y,\pm e_{\ell}) \ {\rm if} \ 0 < y_{\ell}<1.
\end{equation*}
If $y_{\ell}=1$, then $q(k,y,e_{\ell}) =0$ and $q(k,y,-e_{\ell})=p(k,y,e_{\ell})+p(k,y,-e_{\ell})$. If 
$y_{\ell}=0$, then $q(k,y,-e_{\ell}) =0$ and $q(k,y,e_{\ell})=p(k,y,e_{\ell})+p(k,y,-e_{\ell})$.
\vspace{5pt}
\\
The deadlock time of the banker algorithm is then given by $T^{(m)} \equiv \inf
\{ n \geq 0, \sum_{\ell=1}^d \langle R_n,e_{\ell} \rangle \geq \Lambda \}$.
This also writes: $T^{(m)} = \inf \{n \geq 0, \ R_n \in mF_0\}$, with $F_0=\{ y \in [0,2]^d, \  \sum_{i\leq d} |y_i-1| \leq (d-\lambda)\}$ (see Figure
\ref{figure:f0} below). This latter form will be useful in the proof of the main results.
\begin{figure} [htb]
\begin{center}
\includegraphics[
width=0.3\textwidth,angle=0]
{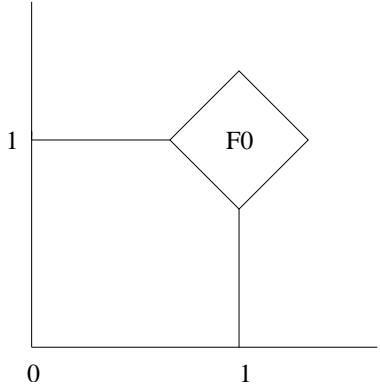}
\caption{Absorption set.}
\label{figure:f0}
\end{center}
\end{figure}
\subsection{Main Assumptions}
In formula \eqref{star}, the 
division by  $m$ indicates that the dependence of the transition kernel
on the position of the walker takes place at scale $m$. For large  $m$,
the space dependence is mild, since we will assume all through the paper
 the following smoothness property:
\vspace{5pt}
\\
{\bf Assumption (A.1).} {\it The function $p$ is twice continuously differentiable with respect to $y$ with bounded derivatives.
In particular, there exists a constant $K >0$ such that:}
\begin{equation*}
\forall x \in E, \ u  \in {\mathcal V}, \ (y,y') \in ({\mathbb R}^d)^2, \ |p(x,y,u) - p(x,y',u)| \leq K |y-y'|.
\end{equation*}
It is then readily seen that the transition kernel (\ref{star}) weakly depends on the space position of the walker:
a step of the walker 
modifies only slightly 
the transition kernel. 
\vspace{5pt}
\\
We also assume the environment to be ergodic and to fulfill the so-called central limit theorem for Markov
chains. We thus impose the following sufficient conditions:
\vspace{5pt}
\\
{\bf Assumption (A.2).} {\it The matrix $P$ is irreducible on $E$. It is then well known that $P$ admits a unique invariant probability 
measure, denoted by $\mu$.}
\vspace{5pt}
\\
{\bf Assumption (A.3)} {\it The matrix $P$ is aperiodic. In particular, it satisfies the Doeblin condition:}
\begin{equation*}
\exists m \geq 0, \ \exists \eta >0, \ \forall (k,\ell) \in E^2, \ (P^m)(k,\ell) \geq \eta.
\end{equation*}
{\bf Assumption (A.4)} {\it Define for $k \in E$ and $y \in \bbR^d$, $g(k,y) \equiv \sum_{u \in {\mathcal V}} [p(k,y,u)u]$ (\emph{i.e.}, $g(k,y)$ matches the
expectation of the measure $p(k,y,\cdot)$) and assume that $g(\cdot,y)$, seen as a function from $E$ to $\bbR^d$, is centered with respect to the measure
$\mu$}.
\vspace{5pt}
\\
Let us briefly describe the role of each of these conditions:
\begin{enumerate}
\item Thanks to Assumption {\bf (A.2)}, the Markov chain $(\xi_n)_{n \geq 0}$ satisfies the ergodic theorem 
for Markov processes. 
\item The general central limit theorem for Markov chains with finite state space follows from the 
 Doeblin condition, given in Assumption {\bf (A.3)}.
\item Assumption {\bf (A.4)} permits to apply the previous central limit theorem to the function $g$.
\end{enumerate}
\subsection{Frequently Used Notations}
For a square integrable martingale $M$, $\langle M \rangle$ denotes the bracket of $M$
(do not mix up with the Euclidean scalar product, which is denoted by $\langle
x,y \rangle$ for $x,y \in {\mathbb R}^d$). 
\vspace{5pt}
\\
We also denote by ${\mathbb D}({\mathbb R}_+,{\mathbb R}^d)$ the path
space of right continuous and left limited functions from ${\mathbb R}_+$ into ${\mathbb R}^d$.
A function in ${\mathbb D}({\mathbb R}_+,{\mathbb R}^d)$ is then said to be ``c\`ad-l\`ag'' for the French acronym ``continue
\`a droite-limite \`a gauche''.
\mysection{Main Results}
\label{Section_Main_Results}
\subsection{Asymptotic Behaviour of the Walk without Reflection}
In light of Assumptions {\bf (A.1-4)}, we expect 
the global effect of the environment process $(\xi_n)_{n \geq 0}$ to
reduce for large time
to a deterministic one. To this end, we view the process $(\bar{S}^{(m)}_t \equiv m^{-1} S_{\lfloor m^2 t \rfloor}^{(m)})_{t \geq 0}$ as a 
random element in the space ${\mathbb D}({\mathbb R}_+,{\mathbb R}^d)$:
\begin{thm}
\label{THM1}
The process $(\bar{S}^{(m)}_t)_{t \geq 0}$ converges in law in ${\mathbb D}({\mathbb R}_+,{\mathbb R}^d)$
endowed with the Skorokhod topology towards the (unique) solution of the martingale problem starting from the origin at time zero and
associated to the operator:
\begin{equation*}
{\mathcal L} \equiv \frac{1}{2} \sum_{i,j=1}^d \bar{a}_{i,j}(y) \frac{\partial^2}{\partial x_i \partial x_j} + \sum_{i=1}^d \bar{b}_i(y) \frac{\partial}{\partial y_i}.
\end{equation*}
The limit coefficients $\bar{a}$ and $\bar{b}$ are given by:
\begin{equation*}
\begin{split}
\forall y \in \bbR^d, \ &\bar{a}(y) \equiv \int_{E} \bigl[ \alpha + g v^t + vg^t - 2gg^t \bigr](i,y) d\mu(i),
\\
&\bar{b}(y) \equiv \int_{E} \bigl[ (\nabla_y v - \nabla_y g)g \bigr](i,y) d\mu(i),
\end{split}
\end{equation*}
where $\nabla$ stands for the gradient and $\alpha(i,y)$ denotes the second order moment matrix of the measure $p(i,y,\cdot)$:
\begin{equation*}
\forall i \in E, \ \forall y \in \bbR^d, \ \alpha(i,y) \equiv \sum_{u \in {\mathcal V}} \bigl[uu^t p(i,y,u) \bigr],
\end{equation*}
and $v(i,y)$ denotes the well-defined sum:
\begin{equation*}
\forall i \in E, \ \forall y \in \bbR^d, \ v(i,y) \equiv \sum_{n \geq 0} \sum_{j \in E} \bigl[ P^n(i,j) g(j,y) \bigr].
\end{equation*}
Moreover, for every $y \in \bbR^d$, the lowest eigenvalue of $\bar{a}(y)$ is greater or equal than the lowest 
eigenvalue of the matrix:
\begin{equation*}
\int_{E} \bigl[ \alpha(i,y) - gg^t(i,y) \bigr] d\mu(i).
\end{equation*}
In particular, if the covariance matrices of the measures $(p(i,y,\cdot))_{i \in  E,y \in \bbR^d}$ are uniformly elliptic, the matrices 
$(\bar{a}(y))_{y \in {\mathbb R}^d}$ are also uniformly elliptic.
\end{thm}
The existence of $v$, $\nabla_y v$ and $\nabla_y g$ will be detailed in the sequel of the paper.
\vspace{5pt}
\\
We can denote by $\bar{\sigma}(x)$, for $x \in \bbR^d$, the nonnegative symmetric square root of $\bar{a}(x)$. 
Then, the SDE (stochastic differential equation) associated to ${\mathcal L}$ writes:
\begin{equation}
\label{SDE}
d\bar{X}_t = \bar{b}(\bar{X}_t) dt + \bar{\sigma}(\bar{X}_t) dB_t,
\end{equation}
where $(B_t)_{t \geq 0}$ denotes a $d$-dimensional Brownian motion on a given filtered probability space $(\Omega,({\mathcal F}_t)_{t \geq 0},{\mathbb P})$.
Under Assumptions {\bf (A.1-4)}, \eqref{SDE} is uniquely strongly solvable.
\vspace{5pt}
\\
Of course, if $p(k,y,\cdot)$ does not depend on $y$ (as in the beginning of Subsection \ref{banker.II}), 
then the drift $\bar{b}$ vanishes in Theorem \ref{THM1} since $b=0$ and the matrix-valued function 
 $\bar{a}$ is constant so that $\bar{X}$ is a non-standard Brownian motion.
\vspace{5pt}
\\
Theorem \ref{THM1} appears actually
as an homogenization property: on the long run, the time inhomogeneous rescaled walk $\bar{S}^{(m)}$ behaves 
like an homogeneous diffusion. Underlying techniques to establish the asymptotic 
behaviour of $(\bar{S}^{(m)}_t)_{t \geq 0}$ are well known in the literature
devoted to this subject (see e.g. Bensoussan et al. \cite{blp 1978} for a review on homogenization in 
periodic structures or Jikov et al. \cite{jko 1994} for a monograph on stochastic homogenization). We will detail a few of them in the sequel of the paper.
\vspace{5pt}
\\
We finally mention that Theorem \ref{THM1} is very close to Theorem 3.1 in Guillotin-Plantard and Schott \cite{guillotin 
schott 2002}. Indeed, any dynamical $(E,{\mathcal A},\mu,T)$, as considered by the previous authors, generates an homogeneous Markov chain
with degenerate transitions: $\forall k \in E, \ P(k,Tk)=1$ (but the state space $E$ is then very large). In this framework, the condition 
$\int_E f_j d\mu = (2d)^{-1}$ in Theorem
3.1 in \cite{guillotin schott 2002} implies that the expectation against $\mu$ of the transitions $p(T^ik,e_j) = f_j(T^ix)$ and $p(T^i_k,-e_j)=
1/d - f_j(T^ix)$ given in Subsection 3.1 in \cite{guillotin schott 2002} vanishes, and thus amounts in some sense to Assumption {\bf (A.4)}
in our paper. In the same way, condition $(H)$ in Theorem 3.1 in \cite{guillotin schott 2002} refers more or less to a ``degenerate'' central
limit theorem for the underlying dynamical system and thus to our standing Assumption {\bf (A.3)}.
\vspace{5pt}
\\
However, our own Theorem \ref{THM1} does not recover exactly the result in Guillotin-Plantard and Schott \cite{guillotin schott 2002}:
due to the degeneracy of the transitions of a dynamical system, condition {\bf (A.3)} cannot be satisfied. Actually, the reader must understand that condition {\bf (A.3)}
is a simple but very strong technical condition to ensure the validity of the so-called ``central limit theorem'' for Markov chains. 
It thus permits to draw up a clear framework in which the stochastic homogenization theory applies, but is obviously far from being optimal (refer to Olla
\cite{olla 2003} for a complete overview on central limit theorems for Markov chains).
\vspace{5pt}
\\
Theorem \ref{THM1} is proved in Section \ref{ProofTHM1}.
\subsection{Asymptotic Behaviour of the Reflected Walk}
The rescaled walk $(\bar{R}^{(m)}_t \equiv m^{-1} R_{\lfloor m^2 t \rfloor}^{(m)})_{t \geq 0}$ 
satisfies the following ``reflected version'' of Theorem \ref{THM1}:
\begin{thm}
\label{THM2}
The process $(\bar{R}^{(m)}_t)_{t \geq 0}$ converges in law in ${\mathbb D}({\mathbb R}_+,{\mathbb R}^d)$
endowed with the Skorokhod topology towards the (unique) solution of the martingale problem with normal reflection  
on the boundary of the hypercube $[0,1]^d$, with zero as initial condition and with ${\mathcal L}$ as underlying operator.
\end{thm}
The reflected SDE associated to the hypercube $[0,1]^d$ and to the operator ${\mathcal L}$ writes:
\begin{equation}
\label{RSDE}
d X_t = \bar{b}(X_t) dt + dH_t - dK_t + \bar{\sigma}(X_t) dB_t.
\end{equation}
We explain now the meaning of the different terms in the r.h.s of \eqref{RSDE}. The quantity $\bar{b}(X_t) dt$ refers to
the drift $\bar{b}(\bar{X}_t)dt$ in \eqref{SDE} and  $\bar{\sigma}(X_t) dB_t$ to the stochastic noise
$\bar{\sigma}(\bar{X}_t) dB_t$ in \eqref{SDE}. The new terms $dH$ and $dK$ stand for the differential elements
of two continuous adapted processes with bounded variation that prevent $X$ to leave the hypercube $[0,1]^d$.
More precisely, $H^{1}$ (the first coordinate of $H$) is an increasing process that pushes
(when necessary) the process $X^{1}$ to the right to keep it above zero. In the same way, 
$K^{1}$ is an increasing process that pushes (again when necessary) the process $X^{1}$ to the left to keep it below one.
For each $\ell \in \{2,\dots,d\}$, $H^{\ell}$ and $K^{\ell}$ act similarly in the direction $e_{\ell}$.
Both processes
$H$ and $K$ act in a minimal way. In particular, there is no reason to push the process $X^{\ell}$ when away from the boundary.
\vspace{5pt}
\\
The least action principle for $H$ and $K$ can be summarized as follows. The process $H^{\ell}$, for $\ell \in \{1,\dots,d\}$, 
does not increase when $X^{\ell}$ is different from zero, and the process $K^{\ell}$, for $\ell \in \{1,\dots,d\}$, does not increase
when $X^{\ell}$ is different from one, {\it i.e.},
\begin{equation}
\label{LEAST}
\int_0^{+ \infty} {\mathbf 1}_{\{X^{\ell}_t >0\}}dH_t^{\ell}=0, \ 
\int_0^{+ \infty} {\mathbf 1}_{\{X^{\ell}_t <1\}}dK_t^{\ell}=0.
\end{equation}
For a complete review on reflected SDE's, we refer the reader to Tanaka \cite{tanaka 1979}, Lions and 
Sznitman \cite{lions sznitman 1984} and Saisho \cite{saisho 1987}.
\vspace{5pt}
\\
The following corollary describes the asymptotic behaviour of the deadlock time and the deadlock point of the algorithm.
We assume the matrix $\alpha-gg^t$ to be uniformly elliptic to ensure $\bar{a}$ to be so. In short, the ellipticity of
$\bar{a}$ avoids any singular behaviour (up to a ${\mathbb P}$-null set) of the trajectories of the limit process $X$.
\begin{cor}
\label{COR1}
Assume that $\lambda \in [1,d)$ and that the matrix $\alpha-gg^t$ is uniformly elliptic. Denote by
$T$ the deadlock time of the limit process $X$:
$T \equiv \inf \{ t \geq 0, \ X_t \in F_0 \}$,
where $F_0$ is given in Subsection \ref{banker.II}.
Then, the sequence of rescaled deadlock times $m^{-2} T^{(m)}$ converges in law towards $T$.
In the same way, the sequence of rescaled deadlock points $m^{-1} R_{T^{(m)}}$ converges in law towards
$X_T$.
\end{cor}
\medskip
Theorem \ref{THM2} and Corollary \ref{COR1} are proved in Section \ref{Proof_THM2}.
\subsection{Behaviour of the Limit Process in Dimension Two}
\begin{figure} [htb]
\begin{center}
\includegraphics[
width=0.3\textwidth,angle=0]
{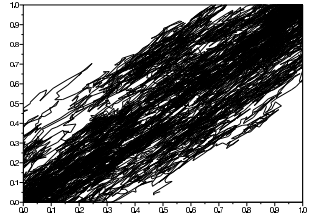}
\caption{One trajectory of the reflected process over $[0,10]$ for $s=0,9$.}
\label{fig:cas1}
\end{center}
\end{figure}

\begin{figure} [htb]
\begin{center}
\includegraphics[
width=0.3\textwidth,angle=0]
{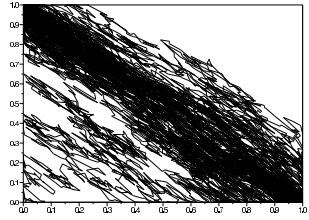}
\caption{One trajectory over $[0,10]$ for $s=-0,9$ (the trajectory is stopped
before $T$).}
\label{fig:cas2}
\end{center}
\end{figure}

\begin{figure} [htb]
\begin{center}
\includegraphics[
width=0.3\textwidth,angle=0]
{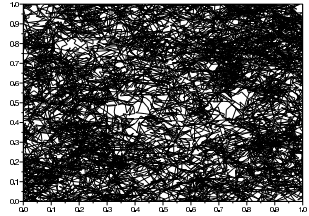}
\caption{One Trajectory over $[0,10]$ for $s=0$.}
\label{fig:cas3}
\end{center}
\end{figure}

A crucial question for numerical applications consists in estimating precisely the mean value of the deadlock time of the limit system.
When the matrix $\overline{a}$ is constant and diagonal and the drift $\bar{b}$ reduces to zero, several explicit computations in terms of Bessel functions are conceivable
(see again \cite{guillotin schott 2002}).
\vspace{15pt}
\\
In our more general framework, the story is quite different. In Delarue \cite{delarue:preprint}, 
 we manage to establish in the two-dimensional
case ({\it i.e.} $d=2$) relevant estimates of the expectation of $T$ 
(now denoted by $T_{\lambda} \equiv \inf \{t \geq 0, \ X_t^1 + X^2_t \geq \lambda \}$
to take into account the parameter $\lambda$), and in particular to distinguish three different asymptotic 
regimes as the parameter $\lambda$ tends to two, each of these regimes depending on the covariance matrix $\bar{a}(1,1)$, and more precisely, on
the sign of its off-diagonal components. 
\vspace{5pt}
\\
Since the matrix $\alpha-gg^t$ is assumed to be uniformly elliptic, the matrix $\bar{a}(1,1)$ is positive and writes:
\begin{equation}
\label{sym_writing}
\bar{a}(1,1)= \left(
\begin{array}{cc}
\rho_1^2 &s\rho_1 \rho_2
\\
s \rho_1 \rho_2 & \rho_2^2
 \end{array}
 \right),
\end{equation}
with $\rho_1,\rho_2>0$ and $s \in ]-1,1[$. The matrix $\bar{a}(1,1)$ admits two eigenvalues:
\begin{equation}
\label{eigenvalues}
\lambda_1 = \frac{1}{2} \bigl[ \rho_1^2 + \rho_2^2 + \delta \bigr], \ 
\lambda_2 = \frac{1}{2} \bigl[ \rho_1^2 + \rho_2^2 - \delta \bigr], \ 
{\rm where} \
\delta \equiv \bigl( \rho_1^4 + \rho_2^4 - 2 (1-2 s^2) \rho_1^2 \rho_2^2 \bigr)^{1/2}.
\end{equation}
Denote by $E_1$ and $E_2$ the associated eigenvectors (up to a multiplicative constant). For $s \not =0$,
\begin{equation*}
E_1 = \left(
\begin{array}{c}
1
\\
(2s \rho_1 \rho_2)^{-1} \bigl( \delta + \rho_2^2 - \rho_1^2 \bigr)
\end{array}
\right),
\
E_2 = \left(
\begin{array}{c}
- (2s \rho_1 \rho_2)^{-1} \bigl( \delta + \rho_2^2 - \rho_1^2 \bigr)
\\
1
\end{array}
\right).
\end{equation*}
Since $\delta+ \rho_2^2 - \rho_1^2 \geq 0$, the signs of the non-trivial coordinates of $E_1$ and $E_2$ are given by the sign of $s$. The
main eigenvector ({\it i.e.} $E_1$) has two positive components for $s>0$, and a positive one and a negative one for $s<0$.
Of course, if $s$ vanishes, $E_1$ and $E_2$ reduce to the vectors of the canonical basis.
\vspace{5pt}
\\
The three different regimes can be distinguished as follows:
\vspace{5pt}
\\
{\bf \emph{Positive Case.}}
\vspace{5pt}
\\
If $s>0$, the main eigenvector of $\bar{a}(1,1)$ ({\it i.e.} $E_1$) is globally oriented from $0$ to the neighbourhood of the corner $(1,1)$ 
and tends to push the limit reflected diffusion towards
the border line. The reflection on the boundary cancels most of the effects of the second eigenvalue and keeps on bringing back the diffusion along the main
axis. As a consequence, the hitting time of the border line is rather small and the following asymptotic holds for the diffusion starting from 0:
\begin{equation*}
\sup_{1 < \lambda < 2} {\mathbb E}(T_{\lambda})< + \infty.
\end{equation*}
This phenomenon is illustrated below (see Figure \ref{fig:cas1}) when $\bar{b}$ reduces to 0 and $\bar{a}$ is the constant matrix given by $\rho_1=\rho_2=1$ and
$s=0,9$. We have plotted there a simulated trajectory of the reflected diffusion process, starting from 0 at time 0, and running from time 0 to time 10 in the box $[0,1]^2$. The algorithm
used to simulate the reflected process is given in Slomi\`nski \cite{slominski 2001}. The eigenvector $E_1$ exactly matches $(1,1)^t$.
\vspace{5pt}
\\
{\bf \emph{Negative Case.}}
\vspace{5pt}
\\
If $s<0$, the main eigenvector of $\bar{a}(1,1)$ is globally oriented from $(1,0)$ to the neighbourhood of the corner $(0,1)$ and attracts the diffusion away
from the border line. Again, the reflection on the boundary cancels most of the effects of the second eigenvalue, and thus, acts now as a trap: the diffusion stays
for a long time along the main axis and hardly touches the boundary.  The hitting time satisfies the following asymptotic behaviour when the diffusion starts from $0$:
\begin{equation*}
\exists c_1,c_2\geq 1, \ \forall \lambda \in ]1,2[, \ c_1^{-1} (2-\lambda)^{-c_2} -c_1 \leq {\mathbb E}(T_{\lambda}) \leq c_1 (2-\lambda)^{-c_2^{-1}} + c_1.
\end{equation*}
This point is illustrated by Figure \ref{fig:cas2} below when $\bar{b}$ vanishes and $\bar{a}$ reduces to the constant matrix $\rho_1=\rho_2=1$
and $s=-0,9$ (again, the initial condition of the process is 0). The eigenvector $E_1$ is given, in this case, by $(1,-1)^t$.
\vspace{5pt}
\\
{\bf \emph{Null Case.}}
\vspace{5pt}
\\
The case $s=0$ is intermediate. Eigenvectors are parallel to the axes and the behaviour of the diffusion is close to the behaviour of the two-dimensional
Brownian motion. For the initial condition $0$:
\begin{equation*}
\exists c_1\geq 1, \ \forall \lambda \in ]1,2[, \ -c_1^{-1} \ln(2-\lambda) - c_1 \leq {\mathbb E}(T_{\lambda}) \leq -c_1 \ln(2-\lambda) + c_1.
\end{equation*}
This is illustrated by Figure \ref{fig:cas3} below when $\bar{b}$ vanishes and $\bar{a}$ reduces to the identity matrix (the initial condition of the process is 0).
\vspace{5pt}
\\
The following theorem sums up these different cases (we refer to Delarue \cite{delarue:preprint} for the whole proof, and just provide a sketch of it
in Section \ref{behaviour}):
\begin{thm}
\label{MAIN_THM2}
Assume that $\alpha- gg^t$ is uniformly elliptic. Then, there exists a constant $C_{\ref{MAIN_THM2}}\geq 1$, depending only on
known parameters and on the ellipticity constant of $\bar{a}$, such that:
\begin{enumerate}
\item For $s>0$, $\displaystyle \sup_{\ell \in ]1,2[} \bbE(T_{\ell}) \leq C_{\ref{MAIN_THM2}}$
\item For $s<0$, set $\displaystyle \beta_- \equiv -s >0, \ \beta_+ \equiv s(s-3)(1+s)^{-1} >0$. Then,
\begin{equation*}
\forall \ell \in ]1,2[, \ C_{\ref{MAIN_THM2}}^{-1} (2-\ell)^{-\beta_-} - C_{\ref{MAIN_THM2}} \leq \bbE(T_{\ell}) \leq C_{\ref{MAIN_THM2}} (2-\ell)^{-\beta_+}
+ C_{\ref{MAIN_THM2}}.
\end{equation*}
\item If $s=0$, $\displaystyle
\forall \ell \in ]1,2[, \ - C_{\ref{MAIN_THM2}}^{-1}  \ln(2-\ell) - C_{\ref{MAIN_THM2}} \leq  
\bbE(T_{\ell}) \leq -C_{\ref{MAIN_THM2}}  \ln(2-\ell)  + C_{\ref{MAIN_THM2}}.$
\end{enumerate}
\end{thm}
Note that Theorem \ref{MAIN_THM2} leaves open many questions. For example, we do not know how to
compute, for $s<0$, the exact value of the ``true'' exponent $\beta \equiv \inf \{ c>0, \ \sup_{\ell \in ]1,2[} \bigl[ (2-\ell)^{c} 
\bbE(T_{\ell}) \bigr] < +
\infty \}$. We are even unable to precise the asymptotic behaviour of $\beta$ as $s \rightarrow -1$ (note indeed
that $\lim_{s \rightarrow -1} \beta_- = 1, \ \lim_{s \rightarrow -1} {\beta}_+ = + \infty$). 
\vspace{5pt}
\\
We have very few ideas about the extension of Theorem
\ref{MAIN_THM2} to the upper dimensional cases. The only accessible case for us is $\overline{a}(1,\dots,1)=I_d$, 
$I_d$ denoting the identity matrix of size $d$: in this case, the analysis derives from the transience properties of the Brownian motion in dimension $d \geq 3$.
\mysection{Proof of Theorem \ref{THM1}}
\label{ProofTHM1}
The strategy of proof of Theorem \ref{THM1} is well-known. 
We aim at writing the sequence $(\bar{S}^{(m)})_{m \geq 1}$ as a sequence of martingales, or at least of semimartingales with relevant drifts. The 
asymptotic behaviour of the martingale parts then follows from a classical central limit argument.
\subsection{Semimartingale Expansion of the Rescaled Walk}
{\it (First Step. A First Martingale Form.)}
\vspace{5pt}
\\
Denote by $({\mathcal F}_n)_{n \geq 0}$ the filtration $({\mathcal F}_n^{\xi,J})_{n \geq 0}$ and write $(S_n)_{n \geq 0}$ as follows:
\begin{equation*}
\begin{split}
S_{n+1} &= S_n + J_{n+1}
\\
&=S_n + {\mathbb E} \bigl(J_{n+1}|{\mathcal F}_n \bigr) + J_{n+1} - {\mathbb E} \bigl(J_{n+1}|{\mathcal F}_n \bigr).
\end{split}
\end{equation*}
The increment $Y_{n+1}^0 \equiv  J_{n+1} -  {\mathbb E}(J_{n+1}| {\mathcal F}_n)$ 
appears as a martingale increment. Referring to \eqref{star} and to the definition of $g$ given in Assumption {\bf (A.4)}, the underlying conditional expectation writes
\begin{equation}
\label{2}
{\mathbb E}  \bigl(J_{n+1}|{\mathcal F}_n \bigr) = g\bigl(\xi_n,S_n/m \bigr).
\end{equation}
Thus, $S_{n+1}$ now writes:
\begin{equation}
\label{star2}
S_{n+1}  = S_n + g\bigl(\xi_n,S_n/m \bigr) + 
Y_{n+1}^0.
\end{equation}
We deduce the following expression for the rescaled walk:
\begin{equation}
\label{star2bis}
\bar{S}^{(m)}_t = m^{-1} \sum_{k=0}^{\lfloor m^2 t \rfloor -1} g \bigl(\xi_k,S_k/m \bigr) + m^{-1} \sum_{k=1}^{\lfloor m^2t \rfloor} 
Y_{k}^0.
\end{equation}
Note at this early stage of the proof that condition {\bf (A.4)} is necessary to establish Theorem \ref{THM1}.
Assume indeed that $g$ does not depend on $S_n/m$. Then, the first term in the above right hand side reduces to
\begin{equation}
\label{star3}
m^{-1} \sum_{k=0}^{\lfloor m^2 t \rfloor -1} g \bigl(\xi_k \bigr).
\end{equation}
The ergodic theorem for Markov chains then yields:
\begin{equation*}
m^{-2} \sum_{k=0}^{\lfloor m^2 t \rfloor -1} g \bigl(\xi_k \bigr) \rightarrow t \int_E g d\mu \ {\rm as} \ m \rightarrow + \infty.
\end{equation*}
In particular, if the expectation of $g$ with respect to $\mu$ doesn't vanish, the term (\ref{star3}) explodes with $m$.
The question is then to study the asymptotic behaviour of \eqref{star3} under the centering condition {\bf (A.4)}, or in other
words to establish a central limit theorem for the sequence $(g(\xi_n))_{n \geq 0}$. 
\vspace{5pt}
\\
{\it (Second Step. Auxiliary Problems.)} 
\label{sec:aux}
\vspace{5pt}
\\
The basic strategy consists in solving the system of Poisson equations driven by the generator of the chain $(\xi_n)_{n \geq 0}$ and by the
$\bR^d$ valued drift $g$. In our framework, this system extends to a family of systems since $g$ does depend on the extra parameter $y \in \bR^d$.
We thus investigate:
\begin{equation}
\label{star4}
(I-P)v=g(\cdot,y), \ y \in {\mathbb R}^d.
\end{equation}
For $y \in \bR^d$, the solvability of this equation is ensured by the Doeblin condition {\bf (A.3)}. Referring to Remark 2, Theorem 4.3.18 in 
Dacunha-Castelle and Duflo \cite{dacunha-castelle duflo t2}, we derive from {\bf (A.3)} that there exist two 
constants $\Gamma$ and $c$, $0 <c<1$, depending only on known
 parameters appearing in Assumptions {\bf (A.1-4)}, such that:
\begin{equation}
\label{star5}
\forall n \geq 0, \ \forall i \in E, \ \sum_{j \in E} |P^n(i,j) - \mu(j)| \leq \Gamma c^n.
\end{equation}
We then claim (see e.g. Theorem 4.3.18 in Dacunha-Castelle and Duflo \cite{dacunha-castelle duflo t2}):
\begin{prop}
\label{PROP1} Under {\bf (A.1-4)}, the equation (\ref{star4}) is solvable for every $y \in {\mathbb R}^d$.
Moreover, there is a unique centered solution with respect to the measure $\mu$. It is given by:
\begin{equation*}
\forall i \in E, \ v(i,y) = \sum_{n \geq 0} \sum_{j \in E} \bigl[P^n(i,j) g(j,y) \bigr].
\end{equation*}
\end{prop}
Note carefully that $v(i,y)$ belongs to $\bR^d$. Note also from {\bf (A.4)} and (\ref{star5}) that $v(i,y)$ is well defined since
\begin{equation*}
 v(i,y) = \sum_{n \geq 0} \sum_{j \in E} \bigl[P^n(i,j)-\mu(j)\bigr] g(j,y) .
\end{equation*}
The question of the regularity of $g$ and $v$ turns out to be crucial in the statement of Theorem \ref{THM1}.
The following Lemmas directly follow from Assumption {\bf (A.1)} and \eqref{star5}:
\begin{lem}
\label{LEMME1}
The function $g(i,y)$ is bounded by 1 and is twice continuously differentiable with respect to $y$ with bounded derivatives. 
In particular, $g$ is Lipschitz continuous in $y$, 
uniformly in $i$. The Lipschitz constant is denoted by $C_{\ref{LEMME1}}$.
\end{lem}
\begin{lem}
\label{LEMME2}
The function $v(i,y)$ is bounded by $\Gamma (1-c)^{-1}$ and is twice continuously differentiable in $y$
with bounded derivatives. In particular, $v$ is Lipschitz continuous in $y$, uniformly in $i$. The Lipschitz constant is denoted by $C_{\ref{LEMME2}}$.
\end{lem}
{\it (Third Step. Final Semimartingale Form.)}
\vspace{5pt}
\\
We claim:
\begin{prop}
\label{PROP1BIS}
Define for $(i,y) \in E \times \bbR^d$:
\begin{equation*}
b^{(m)}(i,y) \equiv \sum_{u \in {\mathcal V}}
\bigl[ p(i,y,u) \bigl( (v-g)(i,y+u/m) - (v-g)(i,y) \bigr)\bigr].
\end{equation*}
and denote by $(Y_n)_{n \geq 1}$ the sequence of corrected martingale increments:
\begin{equation*} 
Y_{n} \equiv J_{n} + v\bigl(\xi_{n},S_{n}/m) - {\mathbb E} \bigl[J_{n} + v\bigl(\xi_{n},S_{n}/m) | {\mathcal F}_{n-1} \bigr].
\end{equation*}
Then, the walk $(S_n)_{n \geq 0}$ writes:
\begin{equation}
\label{star7}
S_{n+1} + v(\xi_{n+1},S_{n+1}/m) = S_n + v(\xi_n,S_n/m)  + b^{(m)}(\xi_n,S_n/m) + Y_{n+1}.
\end{equation}
In particular, the process $(S_n + v(\xi_n,S_n/m))_{n \geq 0}$ is a semimartingale. The term $b^{(m)}$ appears as a drift increment.
It satisfies:
$\forall m \geq 1, i \in E, y \in {\mathbb R}^d, \ |b^{(m)}(i,y)| \leq C_{\ref{PROP1BIS}}/m$,
where $C_{\ref{PROP1BIS}}$ depends only on known parameters in {\bf (A.1-4)}.
\end{prop}
{\bf Proof.} Apply $v$ to the couple $(\xi_{n+1},S_{n+1}/m)$ for a given $n \geq 0$. In this perspective, note that:
\begin{equation*}
\begin{split}
{\mathbb E} \bigl( v(\xi_{n+1},S_{n+1}/m) | {\mathcal F}_n \bigr) &= 
\sum_{i \in E, u \in {\mathcal V}} \bigl[ P(\xi_n,i) p(\xi_n,S_n/m,u)v\bigl(i,(S_n+u)/m \bigr)  \bigr]
\\
&=\sum_{u \in {\mathcal V}} \biggl[ p(\xi_n,S_n/m,u)  \biggl( \sum_{i \in E} P(\xi_n,i) v\bigl(i,(S_n+u)/m \bigr)  \biggr)
\biggr].
\end{split}
\end{equation*}
According to Proposition \ref{PROP1}, we have:
\begin{equation*}
\sum_{i \in E} \bigl[ P(\xi_n,i) v\bigl(i,(S_n+u)/m \bigr) \bigr] = (v-g)(\xi_n,(S_n+u)/m),
\end{equation*}
yielding:
\begin{equation*}
\begin{split}
&{\mathbb E} \bigl( v(\xi_{n+1},S_{n+1}/m) | {\mathcal F}_n \bigr) 
\\
&\hspace{15pt} = 
(v-g)(\xi_n,S_n/m)
+  \sum_{u \in {\mathcal V}} \biggl[  p(\xi_n,S_n/m,u) \biggl( (v-g)\bigl(\xi_n,(S_n+u)/m \bigr) - (v-g)\bigl(\xi_n,S_n/m\bigr) \biggr)\biggr].
\end{split}
\end{equation*}
Due to the definition of $b^{(m)}$, this writes in an equivalent form:
\begin{equation}
\label{star6}
g(\xi_n,S_n/m) = 
v(\xi_n,S_n/m) - {\mathbb E} \bigl( v(\xi_{n+1},S_{n+1}/m) | {\mathcal F}_n \bigr) + b^{(m)}(\xi_n,S_n/m).
\end{equation}
From (\ref{star2}) and (\ref{star6}), derive that:
\begin{equation*}
S_{n+1} = S_n + 
J_{n+1} + v\bigl(\xi_n,S_n/m) - {\mathbb E} \bigl(J_{n+1} + v(\xi_{n+1},S_{n+1}/m) 
|{\mathcal F}_n \bigr)+ b^{(m)}(\xi_n,S_n/m).
\end{equation*}
Thanks to the definition of $(Y_n)_{n \geq 1}$, we recover the required writing. The bound for $b^{(m)}$ derives from
Lemmas \ref{LEMME1} and \ref{LEMME2} (Lipschitz continuity of $g$ and $v$). $\square$
\vspace{5pt}
\\
Investigate the bracket of the martingale part of $(S_n + v(\xi_n,S_n/m))_{n \geq 0}$:
\begin{prop}
\label{PROP1ter}
Denote by $(M_n)_{n \geq 0}$ the square integrable martingale:
$M_n \equiv \sum_{k=1}^n Y_k$.
Then, the bracket of $(M_n)_{n \geq 1}$ is given by:
\begin{equation*}
\langle M \rangle_n = \sum_{k=1}^{n} a^{(m)}(\xi_{k-1},S_{k-1}/m),
\end{equation*}
with (see the statement of Theorem \ref{THM1} for the definition of $\alpha$):
\begin{equation}
\label{star9}
\begin{split}
a^{(m)}(\xi_n,S_n/m) &\equiv \bigl[ \alpha - (v+b^{(m)})(v+b^{(m)})^t \bigr](\xi_n,S_n/m) 
\\
&\hspace{15pt} + \sum_{u \in {\mathcal V}} \biggl[p(\xi_n,S_n/m,u) \bigl[P(vv^t) + u (Pv)^t + (Pv)u^t  \bigr] \bigl(\xi_n,(S_n+u)/m \bigr)  \biggr].
\end{split}
\end{equation}
Moreover, there exists a finite constant $C_{\ref{PROP1ter}}$, depending only on
known parameters in {\bf (A.1-4)}, such that:
$\forall m \geq 1, i \in E, y \in {\mathbb R}^d, \ |a^{(m)}(i,y)| \leq C_{\ref{PROP1ter}}.$
\end{prop}
The notation $(Pv)(i,y)$ in \eqref{star9} (the analogue holds for $(P(vv^t))(i,y)$) stands for the $d$-dimensional vector $(Pv_1(i,y),\dots,Pv_d(i,y))^t$.
\vspace{5pt}
\\
{\bf Proof.} Recall that $|J_n|$ is almost surely bounded by $1$ and that $v$ is a bounded function (see Lemma \ref{LEMME2}). Hence, the variables 
$(Y_n)_{n \geq 1}$ are almost surely bounded and $(M_n)_{n \geq 1}$ is square-integrable.
The bracket is given by:
\begin{equation*}
\langle M \rangle_n = \sum_{k=1}^{n} {\mathbb V} \bigl[  J_{k} + v(\xi_{k},S_ {k}/m)
\bigl| {\mathcal F}_{k-1} \bigr],
\end{equation*}
where ${\mathbb V}[.|{\mathcal F}_{k-1}]$ denotes the conditional covariance with respect to the $\sigma$-field ${\mathcal F}_{k-1}$.
Due to (\ref{2}) and (\ref{star6}), we have:
\begin{equation*}
{\mathbb E} \bigl[ J_{n+1} + v(\xi_{n+1},S_{n+1}/m) |{\mathcal F}_n \bigr] = (v+ b^{(m)})(\xi_n,S_n/m).
\end{equation*}
Since $ {\mathbb V}[ Z |{\mathcal F}_n]={\mathbb E}[ZZ^t|{\mathcal F}_n] - 
{\mathbb E}[Z|{\mathcal F}_n] {\mathbb E}[Z^t|{\mathcal F}_n]$, 
\begin{equation*}
\begin{split}
&{\mathbb V} \bigl[  J_{n+1} + v(\xi_{n+1},S_ {n+1}/m)
\bigl| {\mathcal F}_{n} \bigr]
\\
&=  {\mathbb E} \biggl[ \bigl( J_{n+1} + v(\xi_{n+1},S_ {n+1}/m)  \bigr)  \bigl( J_{n+1} + v(\xi_{n+1},S_ {n+1}/m)   \bigr)^t \bigl| {\mathcal F}_{n} \biggr] 
\\
&\hspace{15pt} - \bigl[ (v+ b^{(m)})(v+ b^{(m)})^t \bigr](\xi_n,S_n/m) 
\\
&=\sum_{i \in E}  \sum_{u \in {\mathcal V}} \biggl[  P(\xi_n,i)p(\xi_n,S_n/m,u) \bigl[u + v \bigl(i,(S_n+u)/m \bigr)  \bigr] \bigl[u + v \bigl(i,(S_n+u)/m \bigr) \bigr]^t \biggr]
\\
&\hspace{15pt} - \bigl[ (v+ b^{(m)})(v+ b^{(m)})^t \bigr](\xi_n,S_n/m).
\end{split}
\end{equation*}
Reducing the sum over $i \in E$:
\begin{equation*}
\begin{split}
{\mathbb V} \bigl[  J_{k} + v(\xi_{k},S_ {k}/m)
\bigl| {\mathcal F}_{k-1} \bigr]
&= \sum_{u \in {\mathcal V}} \bigl[p(\xi_n,S_n/m,u) u u^t  \bigr]
\\
&\hspace{15pt} +  \sum_{u \in {\mathcal V}} \bigl[p(\xi_n,S_n/m,u)  u  (Pv)^t\bigl(\xi_n,(S_n+u)/m\bigr) \bigr] 
\\
&\hspace{15pt} +  \sum_{u \in {\mathcal V}} \biggl[p(\xi_n,S_n/m,u)  (Pv)\bigl(\xi_n,(S_n+u)/m\bigr) u^t  \biggr] 
\\
&\hspace{15pt} + \sum_{u \in {\mathcal V}} \biggl[p(\xi_n,S_n/m,u)  \bigl[P (vv^t)\bigr]\bigl(\xi_n,(S_n+u)/m \bigr)  \biggr]
\\
&\hspace{15pt} - \bigl[ (v+ b^{(m)})(v+ b^{(m)})^t \bigr](\xi_n,S_n/m).
\end{split}
\end{equation*}
This recovers the form of $a^{(m)}$. The bound for $a^{(m)}$ follows from
Lemmas \ref{LEMME1} and \ref{LEMME2}. $\square$
\subsection{Scaling Procedure}
{\it (First Step. Semimartingale Form of the Rescaled Process.)}
\vspace{5pt}
\\
We deduce from Proposition \ref{PROP1BIS} that for all $t \geq 0$ :
\begin{equation*}
\bar{S}^{(m)}_t  + m^{-1} v(\xi_{\lfloor m^2 t\rfloor},\bar{S}_t^{(m)}) 
= m^{-1} v(\xi_{0},0)  + m^{-1} \sum_{k=0}^{\lfloor m^2 t \rfloor -1 } b^{(m)}(\xi_{k},S_{k}/m) + m^{-1} \sum_{k=1}^{\lfloor m^2 t \rfloor} Y_k.
\end{equation*}
Of course, if $\lfloor m^2 t \rfloor =0$, the latter sums reduce to zero.
\vspace{5pt}
\\
Since the function $v$ is bounded (see Lemma \ref{LEMME2}),  the supremum norm of the process $(m^{-1} v(\xi_{\lfloor m^2 t\rfloor},$ 
$\bar{S}_t^{(m)}))_{t \geq 0}$ tends almost surely 
towards 0 as $m$ increases. In particular, $\bar{S}^{(m)}$ writes, up to a negligible term:
\begin{equation} 
\label{187}
\forall t \geq 0, \ \bar{S}^{(m)}_t = m^{-1} \sum_{k=0}^{\lfloor m^2 t \rfloor -1} b^{(m)}(\xi_{k},S_{k}/m) + m^{-1} \sum_{k=1}^{\lfloor m^2 t \rfloor} Y_k + O(1/m).
\end{equation}
Set now for the sake of simplicity:
\begin{equation}
\label{1871}
\bar{B}^{(m)}_t \equiv m^{-1} \sum_{k=0}^{\lfloor m^2 t \rfloor -1} b^{(m)}(\xi_{k},S_{k}/m),
\ \bar{M}^{(m)}_t \equiv m^{-1} \sum_{k=1}^{\lfloor m^2 t \rfloor} Y_k.
\end{equation}
{\it (Second Step. Tightness Properties.)} 
\begin{prop}
\label{PROP1.4}
The family $(\bar{S}^{(m)},\bar{M}^{(m)},\langle \bar{M}^{(m)}\rangle)_{m \geq 1}$ is ${\mathcal C}$-tight in  
${\mathbb D}({\mathbb R}_+,{\mathbb R}^{2d+d^2})$.
\end{prop}
{\bf Proof.} We first establish the tightness of the processes $(\bar{B}^{(m)},\bar{M}^{(m)})_{m \geq 1}$ in the space 
of c\`ad-l\`ag functions endowed with the Skorokhod topology. To this end, note 
from Proposition \ref{PROP1BIS} (estimates for $(b^{(m)})_{m \geq 1}$)
 that for any $(s,t) \in ({\mathbb R}_+)^2$, $0 \leq s <t$:
\begin{equation}
\label{BV}
\begin{split}
|\bar{B}^{(m)}_t - \bar{B}^{(m)}_s| &\leq Cm^{-2} \big( \lfloor m^2t \rfloor - \lfloor m^2 s  \rfloor  \big)
\\
&\leq C \big( (t-s)+m^{-2} \big).
\end{split}
\end{equation}
Therefore, the family $(\bar{B}^{(m)})_{m \geq 1}$ is ${\mathcal C}$-tight in ${\mathbb D}({\mathbb R}_+,{\mathbb R}^d)$ endowed with
the Skorokhod topology: all the limit processes of the sequence have continuous trajectories (see e.g. Jacod and Shiryaev \cite{jacod shiryaev 1987}, Propositions
VI 3.25 and VI 3.26).
\vspace{5pt}
\\
We turn now to the process $\bar{M}^{(m)}$. It is a square integrable martingale with respect to the filtration $({\mathcal F}_{\lfloor m^2 t \rfloor})_{t \geq 0}$.
Moreover, from Proposition \ref{PROP1ter}, the bracket of the martingale is given by:
\begin{equation}
\label{star10}
\langle \bar{M}^{(m)} \rangle_t = m^{-2} \sum_{k=0}^{\lfloor m^2 t \rfloor -1} a^{(m)}(\xi_{k},S_{k}/m).
\end{equation}
Since the functions $(a^{(m)})_{m \geq 1}$ are uniformly bounded, 
the family  $(\langle \bar{M}^{(m)}\rangle)_{m \geq 1}$ is also ${\mathcal C}$-tight
in ${\mathbb D}({\mathbb R}_+,{\mathbb R}^{2d})$, by the same argument as 
above.
\vspace{5pt}
\\
Referring to Theorem VI 4.13 in Jacod and Shiryaev \cite{jacod shiryaev 1987}, we deduce that the family $(\bar{M}^{(m)})_{m \geq 1}$ is tight in  ${\mathbb D}({\mathbb R}_+,{\mathbb R}^{d})$.
It is even ${\mathcal C}$-tight since the jumps of the martingale $\bar{M}^{(m)}$ are bounded in absolute value by $C/m$.
\vspace{5pt}
\\
Finally, we conclude that the family $(\bar{S}^{(m)},\bar{M}^{(m)},\langle \bar{M}^{(m)}\rangle)_{m \geq 1}$ is ${\mathcal C}$-tight in  
${\mathbb D}({\mathbb R}_+,{\mathbb R}^{2d+d^2})$ (see Corollary VI 3.33 in Jacod and Shiryaev \cite{jacod shiryaev 1987}). $\square$
\vspace{5pt}
\\
{\it (Third Step. Extraction of a Converging Subsequence.)}
\begin{lem}
\label{LEMME3}
Every limit $(\bar{X},\bar{M},L)$ (in the sense of the weak convergence in 
${\mathbb D}({\mathbb R}_+,{\mathbb R}^{2d+d^2})$) of the sequence of processes 
$(\bar{S}^{(m)},\bar{M}^{(m)},\langle \bar{M}^{(m)}\rangle)_{m \geq 1}$ is continuous.
Moreover, $\bar{M}$ is a square integrable martingale with respect to the right continuous filtration generated by the process $(\bar{X},\bar{M},L)$
augmented with ${\mathbb P}$-null sets. Its bracket coincides with $L$.
\end{lem}
{\bf Proof.}
According to Proposition IX 1.17 in Jacod and Shiryaev \cite{jacod shiryaev 1987}, $\bar{M}$ is a 
square integrable martingale.  The same argument shows that $\bar{M}^2-L$ is a martingale. It is readily seen from
(\ref{star10}) that $L$ is of bounded variation and satisfies $L_0=0$. Thus, $\langle \bar{M}\rangle=L$ (see Theorem I 4.2 in \cite{jacod shiryaev 1987}). $\square$
\vspace{5pt}
\\
We aim at expressing the process $\langle \bar{M}\rangle$. To this end,
 we have to study the asymptotic behaviour of $\langle \bar{M}^{(m)} \rangle$:
\begin{lem}
\label{LEMME4} Referring to the definition of $\bar{a}$ in the statement of Theorem \ref{THM1}:
\begin{equation*}
 \sup_{0 \leq t \leq T} \biggl[ \biggl| \langle \bar{M}^{(m)} \rangle_t - m^{-2}
 \sum_{k=0}^{\lfloor m^2 t \rfloor -1} \bar{a}(S_k/m) \biggr| \biggr] \overset{{\mathbb P}-{\rm probability}}{\longrightarrow} 0 \ {\rm as}
 \ m \rightarrow + \infty.
\end{equation*}
\end{lem}
{\bf Proof.}
The function $a^{(m)}$ from (\ref{star9}) 
uniformly converges on $E \times {\mathbb R}^d$ towards the function $a$ given by:
\begin{equation}
\label{star12}
\begin{split}
\forall i \in E, \ \forall y \in {\mathbb R}^d, \ a(i,y) 
&= \alpha(i,y)  +  \bigl[ g (Pv)^t + (Pv)g^t + P(vv^t) -vv^t \bigr](i,y)
\\
&=\alpha(i,y) + \bigl[ gv^t+vg^t - 2gg^t + P(vv^t) - vv^t \bigr](i,y).
\end{split}
\end{equation} 
This follows  from the {\it a priori} estimates for $(b^{(m)})_{m \geq 1}$ (see Proposition \ref{PROP1BIS})
and from Lemma \ref{LEMME2}. Thus, $\langle \bar{M}^{(m)} \rangle$ behaves like:
\begin{equation}
\label{star13}
 \biggl(m^{-2} \sum_{k=0}^{\lfloor m^2 t \rfloor -1} a(\xi_k,S_k/m) \biggr)_{t \geq 0}.
\end{equation} 
Lemma \ref{LEMME4} then follows from Proposition \ref{PROP2} below. $\square$
\begin{prop}
\label{PROP2}
Let $f: E \times {\mathbb R}^d \rightarrow {\mathbb R}$ be bounded and uniformly Lipschitz continuous with respect to the second variable.
 Then, we have for all $ T \geq 0$:
\begin{equation*}
 \sup_{0 \leq t \leq T} \biggl[ m^{-2} \biggl|
\sum_{k=0}^{\lfloor m^2 t \rfloor -1} f(\xi_k,S_k/m) - \sum_{k=0}^{\lfloor m^2 t \rfloor -1} \bar{f}(S_k/m) \biggr| \biggr] \overset{{\mathbb P}-{\rm probability}}{\longrightarrow} 0 \ {\rm as}
 \ m \rightarrow + \infty,
\end{equation*}
with $\displaystyle \forall y \in {\mathbb R}^d, \ \bar{f}(y) = \int_E f(i,y) d\mu(i)$.
\end{prop}
The proof of Proposition \ref{PROP2} relies on the ergodic theorem for $(\xi_n)_{n \geq 0}$ and on the tightness of the family $(\bar{S}^{(m)})_{m \geq 1}$. It is postponed to the end of the section.
\vspace{5pt}
\\
{\it (Third Step. Identification of the Limit.)}
\begin{lem}
\label{LEMME5}
Under the assumptions and notations of Lemmas \ref{LEMME3} and \ref{LEMME4} and of Proposition \ref{PROP2}, the bracket of $\bar{M}$ is given by
$\forall t \geq 0, \ \langle \bar{M} \rangle_t = \int_0^t  \bar{a}(\bar{X}_s)ds$,
and the couple $(\bar{X},\bar{M})$ satisfies:
\begin{equation}
\label{star15}
\forall t \geq 0, \ \bar{X}_t = \int_0^t \bar{b}(\bar{X}_s)ds + \bar{M}_t.
\end{equation} 
\end{lem}
{\bf Proof.} Focus first on the limit form of the bracket. Note to this end that for every $t \geq 0$:
\begin{equation}
\label{1872}
m^{-2} \sum_{k=0}^{\lfloor m^2 t \rfloor -1} \bar{a}(S_k/m) = \int_0^t \bar{a}\bigl(\bar{S}^{(m)}_s\bigr)ds + O(1/m).
\end{equation}
The convergence will then follow from the one of $(\bar{S}^{(m)})_{m \geq 0}$.
\vspace{5pt}
\\
In fact, Proposition \ref{PROP2} applies to studying the asymptotic behaviour of the drift part. Derive indeed from Lemmas \ref{LEMME1} and \ref{LEMME2} and 
Proposition \ref{PROP1BIS} that $mb^{(m)}$ uniformly converges on $E \times {\mathbb R}^d$ towards
$(\nabla_y v - \nabla_y g)g$. Thus, the asymptotic behaviour 
of $\bar{B}^{(m)}$ in \eqref{1871} just reduces to the one of:
\begin{equation*}
m^{-2} \sum_{k=0}^{\lfloor m^2 t \rfloor-1} b(\xi_k,S_k/m), \quad 
b \equiv (\nabla_y v - \nabla_y g)g = \bigl( \sum_{j=1}^d \frac{\partial (v^i - g^i)}{\partial y_{j}} g_j \bigr)_{1 \leq i \leq d}.
\end{equation*}
According once again to Proposition \ref{PROP2}, this just amounts to study:
\begin{equation}
\label{186}
\begin{split}
&m^{-2} \sum_{k=0}^{\lfloor m^2 t \rfloor -1} \bar{b}(S_k/m) = \int_0^t \bar{b} \bigl(\bar{S}^{(m)}_s \bigr)ds + O(1/m), 
\\ 
&\bar{b}(y) \equiv  \int_E b(i,y) d\mu(i)=
\int_E \bigl[ (\nabla_y v - \nabla_y g)g \bigr](i,y) d\mu(i).
\end{split} 
\end{equation}
Note now from Lemma \ref{LEMME3} and from the continuity of $\bar{a}$ and $\bar{b}$ (see Lemma \ref{LEMME2} and the proof of Lemma \ref{LEMME4}) that the processes:
\begin{equation}
\label{LEM5.3}
\biggl(\bar{S}^{(m)},\bar{M}^{(m)}, \int_0^. \bar{b} \bigl(\bar{S}^{(m)}_s \bigr)ds,\int_0^.\bar{a} \bigl(\bar{S}^{(m)}_s \bigr)ds \biggr)_{m \geq 1}
\end{equation}
converge in law in ${\mathbb D}({\mathbb R}_+,{\mathbb R}^{3d+d^2})$ towards $\displaystyle
\biggl(\bar{X},\bar{M},\int_0^. \bar{b}(\bar{X}_s)ds,\int_0^. \bar{a}(\bar{X}_s)ds \biggr).$
\vspace{5pt}
\\
From Lemmas \ref{LEMME3} and \ref{LEMME4} and \eqref{1872}, the bracket of $\bar{M}$ is given by:
$\forall t \geq 0, \ \langle \bar{M} \rangle_t = \int_0^t  \bar{a}(\bar{X}_s)ds$,
 and the 
relation (\ref{187}) yields \eqref{star15}. $\square$
\vspace{5pt}
\\
We are now in position to complete the proof of Theorem \ref{THM1}. Thanks to Lemmas \ref{LEMME1} and \ref{LEMME2}, the coefficients
$\bar{a}$ and $\bar{b}$ given in the statement of Theorem \ref{THM1} are once differentiable with bounded derivatives.
It is then well known that $\bar{X}$ in Lemma \ref{LEMME5} turns to be the unique solution (starting from the origin at time zero)
to the martingale problem associated to $(\bar{b},\bar{a})$ (see e.g. Chapter VI in Stroock and Varadhan \cite{stroock varadhan 1979}).
\vspace{5pt}
\\
We now prove that the diffusion matrix $\bar{a}$ is elliptic when $\alpha -gg^t$ is uniformly elliptic.
Referring to the definition of $\bar{a}$ (see Theorem \ref{THM1}), 
it is sufficient to prove that the symmetric matrix:
\begin{equation*}
\int_E \bigl[ vg^t + gv^t - gg^t \bigr](i,y)d\mu(i),
\end{equation*}
is non-negative. Recall that $g=v-Pv$. Thus:
\begin{equation}
\label{16bis}
\begin{split}
vg^t + gv^t - gg^t &= vv^t - v(Pv)^t + vv^t - (Pv)v^t - vv^t - (Pv)(Pv)^t + v(Pv)^t + (Pv)v^t
\\
&= vv^t - (Pv)(Pv)^t.
\end{split}
\end{equation}
We see from Schwarz inequality 
that for every $(x,y) \in {\mathbb R}^d \times {\mathbb R}^d$:
\begin{equation}
\label{16ter}
\begin{split}
\int_E \langle x,(Pv)(Pv)^t x \rangle(i,y) d\mu(i)&= \int_E \langle x,Pv \rangle^2(i,y) d\mu(i)
\\
&= \int_E \biggl( \sum_{j \in E} P(i,j)\langle x,v(j,y)\rangle \biggr)^2 d\mu(i)
\\
&\leq \int_E \sum_{j \in E} P(i,j) \langle x,v(j,y) \rangle^2 d\mu(i)
\\
&= \int_E \langle x,v(i,y) \rangle^2 d\mu(i)
\\
&= \int_E \langle x,(vv^t) x \rangle (i,y) d\mu(i).
\end{split}
\end{equation}
Thanks to (\ref{16bis}) and (\ref{16ter}), we complete the proof of Theorem \ref{THM1}. $\square$\vspace{5pt}
\subsection{Proof of Proposition \ref{PROP2}} 
The proof follows the strategy in Pardoux and Veretennikov \cite{pardoux veretennikov 1997}. For this reason, we just present a sketch of it.
\vspace{5pt}
\\
Let $T$ be an arbitrary positive real number. 
Recall that $(\bar{S}^{(m)})_{m \geq 1}$ is tight in ${\mathbb D}({\mathbb R}_+,{\mathbb R}^d)$. 
Hence, for a given $\varepsilon>0$, there exists a compact set $A \subset {\mathbb D}({\mathbb R}_+,{\mathbb R}^d)$ such that,
for every $m \geq 1$, $\bar{S}^{(m)}$ belongs to $A$ up to the probability $\varepsilon$:
\begin{equation}
\label{17bis}
{\mathbb P}\{ \bar{S}^{(m)} \in A \} \geq 1-\varepsilon.
\end{equation}
Since $A$ is compact for the Skorokhod topology, we know, for $\delta$ small enough, that every function $\ell \in A$ satisfies:
\begin{equation}
\label{star17}
w_{T+1}'(\ell,\delta) < \varepsilon,
\end{equation}
where $w_{T+1}'$ denotes the usual modulus of continuity on $[0,T+1]$ of a c\`ad-l\`ag function (see e.g. Jacod and Shiryaev, VI 1.8 and Theorem VI 1.14).
Moreover, there exist $\ell_1, \dots, \ell_{N_0}$ 
in $A$ such that, for every $\ell \in A$,
we can find an integer $i \in \{1,\dots,N_0\}$ and a strictly increasing
continuous  function $\lambda_i$ from $[0,T]$ into $[0,T+1/2]$ such that:
\begin{equation}
\label{1873}
\sup_{t \in [0,T]} \bigl[|t - \lambda_i(t)| \bigr] < \delta^2, \ \sup_{t \in [0,T]} \bigl[ \bigl|\ell(t) - \ell_i\bigl(\lambda_i(t)\bigr) \bigr| \bigr] < \delta^2.
\end{equation}
Since $\ell_i$ satisfies (\ref{star17}), there exists a step function $y_i : [0,T+1/2] \rightarrow \bbR$, with all steps but maybe the last one of length greater 
than $\delta$, such that the supremum norm between $y_i$ and $\ell_i$ is less than $\varepsilon$ on $[0,T+1/2]$. The second term in \eqref{1873} then yields: 
$\sup_{t \in [0,T]} [|\ell(t) - y_i(\lambda_i(t))|] < \delta^2 + \varepsilon$.
\vspace{5pt}
\\
Choose now $\omega \in \Omega$ such that $\bar{S}^{(m)}(\omega) \in A$.
Deduce that there exists $i \in \{1,\dots,N_0\}$ such that the supremum on $[0,T]$ of the distance between 
$\bar{S}^{(m)}$ and $y_i \circ \lambda_i$ is less than $\varepsilon + \delta^2$. In particular:
\begin{equation*}
\sup_{k \in [0,\lfloor m^2T-1\rfloor]} \bigl| m^{-1} S_{k} - y_i\bigl(\lambda_i(k/m^2)\bigr) \bigr| \leq 
\sup_{t \in [0,T]} \bigl| \bar{S}^{(m)}_t - y_i\bigl(\lambda_i(t)\bigr) \bigr| \leq \varepsilon + \delta^2.
\end{equation*}
Since $f$ is Lipschitz continuous,
\begin{equation}
\label{18}
\sup_{0 \leq t \leq T} \biggl[ m^{-2}  \sum_{k=0}^{\lfloor m^2 t \rfloor-1}  \bigl|
f(\xi_k,S_k/m) -  f\bigl(\xi_k,y_i\bigl(\lambda_i(k/m^2)\bigr)\bigr) \bigr| \biggr] \leq CT 
\bigl(\varepsilon+\delta^2 \bigr).
\end{equation}
Define now by $t_0=0, \ t_1, \dots, t_p =T$ a subdivision of $[0,T]$ with respect to the steps of $y_i$. Recall that $t_{j}-t_{j-1} \geq \delta$ for $j \in \{1,\dots,p-1\}$.
According to \eqref{1873}, for $t \in G =  \cup_{k=1}^{p-1} [t_{k-1}+\delta^2,t_k-\delta^2] \cup ]t_{p-1}+\delta^2 \wedge T,T]$, $\lambda_i(t)$ and $t$ belong to
the same class of the subdivision, so that $y_i(\lambda_i(t))=y_i(t)$.
Hence,
\begin{equation}
\label{19}
\begin{split}
&\limsup_{m \rightarrow + \infty} \biggl[ m^{-2}  \sum_{k=0}^{\lfloor m^2 T \rfloor-1} \bigl| f\bigl(\xi_k,y_i\bigl(\lambda_i(k/m^2) \bigr)\bigr) 
- f\bigl(\xi_k,y_i(k/m^2) \bigr) \bigr| \biggr]
\\
&\hspace{15pt} \leq C \limsup_{m \rightarrow + \infty} \biggl[ m^{-2} \sum_{k=0}^{\lfloor m^2T \rfloor -1} {\mathbf 1}_{[0,T] \setminus G}(k/m^2) \biggr]
\\
&\hspace{15pt} = C (T-\mu(G)) \leq 2C \delta^2 p  \leq 2C T \delta.
\end{split}
\end{equation}
It follows from (\ref{17bis}), (\ref{18}) and (\ref{19}) that there exists a constant $C$ (not depending on $\delta$, $\varepsilon$, $m$ and $N_0$) such that for $m$ large enough:
\begin{equation}
\label{20}
{\mathbb P} \biggl( \bigcup_{i \in \{1,\dots,N_0\}} \biggl\{ \sup_{0 \leq t \leq T} \biggl[ m^{-2} 
\sum_{k=0}^{\lfloor m^2 t \rfloor-1} \bigl| f(\xi_k,S_k/m) - f\bigl(\xi_k,y_i(k/m^2)\bigr) \bigr|\biggr]
\geq C(\delta + \varepsilon)
\biggr\} \biggr) \leq \varepsilon.
\end{equation}
Note that the same estimate 
holds with $|f(\xi_k,S_k/m)-f(\xi_k,y_i(k/m^2))|$ replaced by $|f(\xi_k,S_k/m)-f(\xi_k,y_i(k/m^2))|+|\bar{f}(S_k/m)-\bar{f}(y_i(k/m^2))|$.
Moreover, for all $i \in E$:
 \begin{equation*}
 \begin{split}
\sum_{k=1}^{\lfloor m^2 T \rfloor} f\bigl(\xi_k,y_i(k/m^2) \bigr) 
&= \sum_{j=1}^p \sum_{k=\lfloor m^2t_{j-1} \rfloor+1}^{\lfloor m^2 t_j \rfloor} f(\xi_k,y_i(k/m^2))
 \\
 &= \sum_{j=1}^p \sum_{k=\lfloor m^2t_{j-1} \rfloor+1}^{\lfloor m^2 t_j \rfloor} f(\xi_k,y_i(t_{j-1})).
\end{split}
\end{equation*}
The ergodic theorem for $(\xi_n)_{n \geq 0}$ yields:
\begin{equation*}
m^{-2}  \sum_{k=1}^{\lfloor m^2 T \rfloor} f(\xi_k,y_i(k/m^2)) \overset{{\mathbb P} \ 
{\rm a.s.}}{\longrightarrow} \sum_{j=1}^p \bigl[ (t_j-t_{j-1}) \bar f(y_i(t_{j-1})) \bigr]
= \int_0^T \bar{f}(y_i(s))ds.
\end{equation*}
Obviously, $m^{-2} \sum_{k=0}^{\lfloor m^2 T \rfloor -1} \bar{f}(y_i(k/m^2) \rightarrow \int_0^T \bar{f}(y_i(s))ds$ as $m$ tends to $+ \infty$.
Hence, for $m$ large enough:
\begin{equation*}
{\mathbb P} \biggl\{ \biggl| 
m^{-2}  \sum_{k=0}^{\lfloor m^2 T \rfloor-1} f(\xi_k,y_i(k/m^2)) - m^{-2} \sum_{k=0}^{\lfloor m^2 T \rfloor-1} \bar{f}(y_i(k/m^2))
\biggr| \geq \varepsilon \biggr\} \leq N_0^{-1} \varepsilon^2.
\end{equation*}
This actually holds for every $t \in [0,T]$. Choosing $t$ of the form $t= k\varepsilon$, $k \in \{0,\dots,\lfloor \varepsilon^{-1}T \rfloor\}$, it comes (for
$m$ large):
\begin{equation*}
{\mathbb P} \biggl\{ \sup_{t \in \{0,\dots,\varepsilon \lfloor \varepsilon^{-1} T \rfloor\}} \biggl| 
m^{-2}  \sum_{k=0}^{\lfloor m^2 t \rfloor-1} f(\xi_k,y_i(k/m^2)) - m^{-2} \sum_{k=0}^{\lfloor m^2 t \rfloor-1} \bar{f}(y_i(k/m^2))
\biggr| \geq \varepsilon \biggr\} \leq  T N_0^{-1} \varepsilon.
\end{equation*}
Since $f$ is bounded, it comes finally for a suitable constant $C>0$ (and for $m$ large):
\begin{equation}
\label{22}
{\mathbb P} \biggl\{ \sup_{0 \leq t \leq T} \biggl| 
m^{-2}  \sum_{k=0}^{\lfloor m^2 t \rfloor-1} f(\xi_k,y_i(k/m^2)) - m^{-2} \sum_{k=0}^{\lfloor m^2 t \rfloor-1} \bar{f}(y_i(k/m^2))
\biggr| \geq C \varepsilon \biggr\} \leq  T N_0^{-1} \varepsilon,
\end{equation}
The sum over $i \in \{1,\dots,N_0\}$ in \eqref{22} is bounded by $\varepsilon$.
Combining (\ref{20}) and (\ref{22}), we complete the proof. $\square$
\mysection{Proofs of Theorem \ref{THM2} and Corollary \ref{COR1}}
\label{Proof_THM2}
\subsection{Proof of Theorem \ref{THM2}}
The proof of Theorem \ref{THM2} is rather similar to the one of Theorem \ref{THM1}. The main difference comes from the
definition of the mean function $g$. In the new setting, it writes for $(i,y) \in E \times [0,1]^d$ $h(i,y) \equiv \sum_{u \in {\mathcal V}}
u q(i,y,u)$. In particular, for $\ell \in \{1,\dots,d\}$:
\begin{equation}
\label{REF0}
\begin{array}{l}
h_{\ell}(i,y)=g_{\ell}(i,y), \ {\rm for} \ 0 < y_{\ell}<1, 
\\
h_{\ell}(i,y) = p(i,y,e_{\ell}) + p(i,y,-e_{\ell}), \ {\rm for} \ y_{\ell}=0,
\\
h_{\ell}(i,y)= - (p(i,y,e_{\ell})+p(i,y,-e_{\ell})), \ {\rm for} \ y_{\ell}=1,
\end{array}
\end{equation}
where $g_{\ell}$ (resp. $h_{\ell}$) denotes the $\ell$th coordinate of $g$ (resp. $h$). 
The process $R$ satisfies the analogue of \eqref{star2bis} with $g$ replaced by $h$. Due to \eqref{REF0},
the drift $h_{\ell}$, for $\ell \in \{1,\dots,d\}$, pushes upwards the $\ell$th coordinate of $R$ when matching 0, and downwards when
matching 1. 
\vspace{5pt}
\\
The drifts of the processes $R$ and $S$ just differ in the extra term $(h-g)$. The following lemma derives from \eqref{REF0}:
\begin{lem}
\label{LEMME1BISII}
For $\ell \in \{1,\dots,d\}$ and $y \in [0,1]^d$, $y_{\ell} = 0 \Rightarrow (h-g)_{\ell} \geq 0$,  
$y_{\ell} = 1 \Rightarrow (h-g)_{\ell} \leq 0$ and $y_{\ell} \in ]0,1[ \Rightarrow (h-g)_{\ell} = 0$.
\end{lem}
(\emph{First Step. Semimartingale Form.})
\vspace{5pt}
\\
Of course, the function $h$ does not satisfy the centering condition {\bf (A.4)}. However, for $g$ as in Section \ref{ProofTHM1},
we can still consider the solution $v$ to the family of Poisson equations \eqref{star4}. Following the proof of Proposition \ref{PROP1BIS}, 
it is plain to derive:
\begin{prop}
\label{PROP1BISII}
Define $c^{(m)}$ as $b^{(m)}$ in Proposition \ref{PROP1BIS}, but with $p(i,y,u)$ replaced by $q(i,y,u)$. Then, the walk $(R_n)_{n \geq 0}$ satisfies:
\begin{equation}
\label{REF1}
\begin{split}
R_{n+1} + v(\xi_{n+1},R_{n+1}/m) &= R_n + v(\xi_n,R_n/m) + (h-g)(\xi_n,R_n/m) 
\\
&\hspace{15pt} + c^{(m)}(\xi_n,R_n/m) + Z_{n+1},
\end{split}
\end{equation}
where $\bbE(Z_{n+1}|\cF_{n}^{\xi,R})=0$. The drift increment $c^{(m)}$ satisfies the same bound as $b^{(m)}$: $\forall m \geq 1,i \in E,y \in \bR^d, \ |c^{(m)}(i,y)| \leq C_{\ref{PROP1BIS}}/m$. 
\end{prop}
We are then able to control the variation of $R$:
\begin{lem}
\label{LEMME2BISII}
Define ${\mathbf e}=e_1 + \dots + e_d$. Then, there exists a constant $C_{\ref{LEMME2BISII}}$, depending only on known parameters in {\bf (A.1-4)}, such that 
for $1 \leq n \leq p \leq q$,
\begin{equation}
\label{REF1B}
\begin{split}
&|R_q-R_n|^2 + 2 \sum_{k=p}^{q-1} \langle m {\mathbf e} - R_n, (h-g)^-(\xi_k,R_k/m) \rangle + 2 \sum_{k=p}^{q-1}
\langle R_n,(h-g)^+(\xi_k,R_k/m) \rangle
\\
&\hspace{15pt} \leq |R_p-R_n|^2 + C_{\ref{LEMME2BISII}}(q-p+m) + 2 \sum_{k=p}^{q-1} \langle R_k-R_n,Z_{k+1} \rangle.
\end{split}
\end{equation}
\end{lem}
{\bf Proof.} For $k \in \{n,\dots,q-1\}$, derive from \eqref{REF1}:
\begin{equation}
\label{REF3BIS-1}
\begin{split}
|R_{k+1}-R_n|^2 &= |R_k-R_n|^2 + 2 \langle R_{k+1}-R_k,R_k-R_n \rangle + |R_{k+1}-R_k|^2
\\
&= |R_k-R_n|^2 + 2 \langle R_k - R_n,(h-g)(\xi_k,R_k/m) \rangle
\\
&\hspace{15pt} - 2 \langle R_k - R_n,v(\xi_{k+1},R_{k+1}/m) - v(\xi_k,R_k/m) \rangle 
\\
&\hspace{15pt} + 2 \langle R_k -R_n,c^{(m)}(\xi_k,R_k/m) \rangle + 2 \langle R_k - R_n,Z_{k+1} \rangle
\\
&\hspace{15pt} + |R_{k+1}-R_k|^2.
\\
&\equiv |R_k-R_n|^2 + 2T(1,k) + 2 T(2,k) + 2 T(3,k) + 2 T(4,k) + T(5,k).
\end{split}
\end{equation}
Apply Lemma \ref{LEMME1BISII}, and derive that $- T(1,k)$ writes $\langle m {\mathbf e} - R_n, (h-g)^-(\xi_k,R_k/m) \rangle$ + $\langle R_n,(h-g)^+(\xi_k,R_k/m) \rangle$.
Summing $T(1,k)$, for $k$ running from $p$ to $q-1$, we obtain the second and third terms in the l.h.s of \eqref{REF1B}.
\vspace{5pt}
\\
Perform now an Abel transform to handle $\sum_{k=p}^{q-1} T(2,k)$. The boundary conditions are bounded by $C m$ since $R$ takes its values in $[0,m]^d$ and $v$ is bounded
(see Lemma \ref{LEMME2}). Moreover, the sum $\sum_{k=p+1}^{q-1} \langle R_k-R_{k-1},v(\xi_k,R_k/m) \rangle$ is bounded by $C(q-p)$ since the jumps of $R$ are bounded by 1. Hence,
the term $\sum_{k=p}^{q-1} T(2,k)$ is bounded by $C(q-p+m)$.
\vspace{5pt}
\\
Due to the bound for $c^{(m)}$ in Proposition \ref{PROP1BISII}, the sum $\sum_{k=p}^{q-1} T(3,k)$ is bounded by $C(q-p)$. The sum of $T(4,k)$ over $k$ provides
the third term in the r.h.s of \eqref{REF1B}. 
\vspace{5pt}
\\
Finally, since the jumps of $R$ are bounded by 1, the sum of $T(5,k)$ over $k$ is bounded by $(q-p)$. $\square$
\vspace{5pt}
\\
We turn now to the martingale part in \eqref{REF1}. Following the proof of Proposition \ref{PROP1ter}, we claim:
\begin{prop}
\label{PROPREF2}
Define the martingale $N$ by $\forall n \geq 1, \ N_n \equiv \sum_{k=1}^n Z_k$. Its bracket $\langle N \rangle_n$ satisfies a suitable version of
Proposition \ref{PROP1ter} with $a^{(m)}$ replaced by $\tilde{a}^{(m)}$, given by:
\begin{equation}
\label{REF3}
\begin{split}
\tilde{a}^{(m)}(\xi_n,R_n/m) &= \bigl[ \alpha - (v+h-g+c^{(m)})(v+h-g+c^{(m)})^t \bigr](\xi_n,R_n/m)
\\
&\hspace{15pt} + \sum_{u \in {\mathcal V}} \biggl[ q(\xi_n,S_n/m,u) \bigl[P(vv^t) + u(Pv)^t + (Pv)u^t \bigr](\xi_n,(R_n+u)/m) \biggr].
\end{split}
\end{equation}
\end{prop}
(\emph{Second Step. Tightness Properties.})
\vspace{5pt}
\\
The scaling procedure then applies as follows. Define as in \eqref{1871}:
\begin{equation*}
\begin{split}
&\bar{C}_t^{(m)} = m^{-1} \sum_{k=0}^{\lfloor m^2 t \rfloor -1} c^{(m)}(\xi_k,R_k/m), \ \bar{N}_t^{(m)} = m^{-1} N_{\lfloor m^2 t \rfloor}, 
\\
&\bar{H}_t^{(m)} = m^{-1} \sum_{k=0}^{\lfloor m^2 t\rfloor-1} (h-g)^+(\xi_k,R_k/m), \ \bar{K}_t^{(m)} 
= m^{-1} \sum_{k=0}^{\lfloor m^2 t \rfloor-1} (h-g)^-(\xi_k,R_k/m),
\end{split}
\end{equation*}
\begin{prop}
The family $(\bar{R}^{(m)},\bar{N}^{(m)},\langle \bar{N}^{(m)} \rangle,\bar{H}^{(m)},\bar{K}^{(m)})_{m \geq 1}$ is $\cC$-tight 
in ${\mathbb D}(\bbR_+,\bbR^{4d+d^2})$.
\end{prop}
{\bf Proof.}
Following Proposition \ref{PROP1.4}, the family $(\bar{N}^{(m)},\langle \bar{N}^{(m)} \rangle)_{m \geq 1}$ is $\cC$-tight in 
${\mathbb D}(\bbR_+,\bbR^{d+d^2})$. 
\vspace{5pt}
\\
Apply now Lemma \ref{LEMME2BISII} to establish the same property for $(\bar{R}^{(m)})_{m \geq 1}$. Choose  
$n=p=\lfloor m^2s \rfloor$ and $q=\lfloor m^2t\rfloor$, with $s<t$. Since all the terms in the l.h.s of \eqref{REF1B} are 
nonnegative, it comes:
\begin{equation}
\label{REF1C}
| \bar{R}^{(m)}_{t}-\bar{R}^{(m)}_s |^2   \leq C (t-s) + 2 m^{-2} \sum_{k=\lfloor m^2 s \rfloor}^{\lfloor m^2 t \rfloor -1}  \langle R_k, Z_{k+1} \rangle + O(1/m).
\end{equation}
The second term in the r.h.s of the above inequality defines a family of martingales. It is plain to see that it is ${\mathcal C}$-tight in ${\mathbb D}(\bR_+,\bR)$ and thus to
derive the ${\mathcal C}$-tightness of the family $(\bar{R}^{(m)})_{m \geq 1}$ in ${\mathbb D}(\bR_+,\bR^d)$. 
\vspace{5pt}
\\
We turn finally to the tightness of $(\bar{H}^{(m)},\bar{K}^{(m)})_{m \geq 1}$. Choose $p$ and $q$ as above and $n=0$ in \eqref{REF1B}, i.e. $R_n=0$. It comes:
\begin{equation*}
\begin{split}
\langle {\mathbf e},\bar{K}^{(m)}_t - \bar{K}^{(m)}_s \rangle &\leq -(\bar{R}_{t}^{(m)})^2 + (\bar{R}_{s}^{(m)})^2 
\\
&\hspace{15pt} + C(t-s+1/m) + 2 m^{-2} \sum_{k=\lfloor m^2s \rfloor}^{\lfloor m^2t \rfloor-1}  \langle R_k,Z_{k+1} \rangle.
\end{split}
\end{equation*}
Note that $|(\bar{R}_{t}^{(m)})^2-(\bar{R}_{s}^{(m)})^2| \leq 2 |\bar{R}_{t}^{(m)}-
\bar{R}_{s}^{(m)}|$ since $\bar{R}^{(m)}$ belongs to $[0,1]^d$. Due to \eqref{REF1C}, deduce that $(\bar{K}^{(m)})_{m \geq
1}$ is ${\mathcal C}$-tight in ${\mathbb D}(\bR_+, \bR^d)$. Of course, the same holds for $(\bar{H}^{(m)})_{m \geq 1}$.
\vspace{5pt}
\\
(\emph{Third Step. Extraction of a Converging Subsequence and Identification of the Limit.})
\begin{prop}
\label{PROPREF2B}
Any weak limit $(X,\tilde{N},H,K)$ of the sequence $(\bar{R}^{(m)},\bar{N}^{(m)},\bar{H}^{(m)},\bar{K}^{(m)})_{m \geq 1}$ satisfies:
\begin{equation}
\label{REF1D}
\forall t \geq 0, \ X_t = \int_0^t \bar{b}(X_s)ds + H_t - K_t + \tilde{N}_t,
\end{equation}
where $\tilde{N}$ is a square integrable continuous martingale whose bracket writes $\langle \tilde{N} \rangle_t = \int_0^t \bar{a}(X_s)ds$, $H$
and $K$ are two nondecreasing continuous processes, matching 0 at zero and satisfying condition \eqref{LEAST}.
\end{prop}
{\bf Proof.} Note first that Proposition \ref{PROP2} still applies in the reflected setting, with $S$ replaced by $R$.
\vspace{5pt}
\\
Focus now on the asymptotic behaviour of the drift $c^{(m)}$ in \eqref{REF1}. As $mb^{(m)}$ does, $mc^{(m)}$ uniformly converges towards 
$(\nabla_y v - \nabla_y g)g$. Following \eqref{186}, this provides the form of the limit drift.
\vspace{5pt}
\\
Concerning the martingale part, we now prove the analogue of Lemma \ref{LEMME4}.  The point is to study the asymptotic form
of the bracket $\langle \bar{N}^{(m)} \rangle$. It well seen from \eqref{REF3} that $\tilde{a}^{(m)}$ uniformly converges towards $a$ (as in \eqref{star12}) 
plus a corrector term that writes $(g-h)\gamma^t + \gamma(g-h)^t$ for a suitable bounded function $\gamma$. Of course, the part in $a$ satisfies
Proposition \ref{PROP2}. The second part vanishes since:
\begin{equation*}
m^{-2} \sum_{k=0}^{\lfloor m^2 t \rfloor -1} |h-g|(\xi_k,R_k/m) \leq m^{-1} \bigl( \bar{H}^{(m)}_t + \bar{K}^{(m)}_t \bigr).
\end{equation*}
This provides the limit form of the brackets $(\langle \bar{N}^{(m)} \rangle)_{m \geq 1}$ and thus the form of the limit bracket.
\vspace{5pt}
\\
We turn finally to the processes $H$ and $K$. It is clear that $H_0=K_0=0$ and that $H$ and $K$ are nondecreasing. It thus remains to 
verify \eqref{LEAST}. Focus on the case of $H$.
Since $(\bar{R}^{(m)},\bar{H}^{(m)})_{m \geq 1}$ converges weakly (up to a subsequence) in the Polish space ${\mathbb D}(\bR,\bR^{2d})$ towards $(X,H)$, we can assume
without loss of generality from the Skorokhod representation theorem (see Theorem II.86.1 in Rogers and Williams \cite{rogers williams T1}) that the convergence holds
almost surely. In particular, for $T>0$, for $\ell \in \{1,\dots,d\}$, for almost every $\omega \in \Omega$, $\langle \bar{H}^{(m)},e_{\ell}\rangle$, seen as a distribution over $[0,T]$ 
converges weakly towards $\langle H,e_{\ell} \rangle$ 
seen as a distribution over $[0,T]$. Hence, for a continuous function $f$ on $[0,T]$:
\begin{equation}
\label{REF4-1}
\lim_{n \rightarrow + \infty} \int_0^T f(t)d \langle \bar{H}^{(m)},e_{\ell} \rangle_t = \int_0^T f(t) dH^{\ell}_t.
\end{equation}
Since $X$ is continuous, \eqref{REF4-1} holds with $f(X_t)$ instead of $f(t)$.
Recall moreover from Proposition VI.1.17 in Jacod and Shiryaev \cite{jacod shiryaev 1987} that $\sup_{0 \leq t \leq T}|\bar{R}^{(m)}_t - X_t|$ tends to 0
almost surely. Since the supremum $\sup_{m \geq 1} \bar{H}_T^{(m)}$ is almost surely finite, \eqref{REF4-1} yields:
\begin{equation}
\label{REF4}
\lim_{n \rightarrow + \infty} \int_0^T f(\bar{R}^{m}_t)d \langle \bar{H}^{(m)},e_{\ell}\rangle_t = \int_0^T f(X_t) d H^{\ell}_t.
\end{equation}
From Lemma \ref{LEMME1BISII}, for a continuous function $f$ with support included in $]0,1]$, the term
$\int_0^T f(X_t^{\ell})d H^{\ell}_t$ vanishes. It is plain to derive
the first equality in \eqref{LEAST}. $\square$
\vspace{5pt}
\\
According to Proposition \ref{PROPREF2}, $(X,H,K)$ satisfies the martingale problem associated to $(\bar{b},\bar{a})$ with normal reflection on 
$\partial [0,1]^d$. Thanks to Theorem 3.1 in
Lions and Sznitman \cite{lions sznitman 1984}, this martingale problem is uniquely solvable. This completes the proof of Theorem \ref{THM2}. $\square$
\subsection{Proof of Corollary \ref{COR1}}
The proof of Corollary \ref{COR1} relies on the mapping Theorem 2.7, Chapter I, in Billingsley \cite{billingsley 1968}. Indeed, consider
the mapping $\Psi: x \in {\mathbb D}(\bR_+,\bR^d) \mapsto \inf \{t \geq 0, x_t \in F_0\}$ (inf $\emptyset$ = $+\infty$). If $x$ denotes a continuous function from $\bR_+$ into
$\bR^d$ and $(x_n)_{n \geq 1}$ a sequence of c\`ad-l\`ag functions from $\bR_+$ into $\bR^d$ converging towards $x$ for the Skorokhod topology,
we know from Proposition VI.1.17 in Jacod and Shiryaev \cite{jacod shiryaev 1987} that $(x_n)_{n \geq 1}$ converges uniformly towards $x$ on compact subsets of $\bR_+$. 
Assume now that
$\Psi(x)$ is finite and that for every $\eta >0$ there exists $t \in ]\Psi(x),\Psi(x)+\eta[$ such that $x(t)$ belongs to the interior of $F_0$ (the point $x(\Psi(x))$ is then said to be 
regular).
Then, $\Psi(x_n)$ tends to $\Psi(x)$. In particular, if almost every trajectory of $X$ satisfies these conditions, $\Psi(\bar{R}^{(m)})$ 
converges in law towards $\Psi(X)$.
\vspace{5pt}
\\
To prove that almost every trajectory of $X$ satisfies the latter conditions, consider the solution $\hat{X}$ 
to the martingale problem associated to the couple $(\hat{b},\hat{a})$, where $\hat{b}$ and $\hat{a}$ are periodic functions of period two in 
each direction of the space given by $(\hat{b}(y),\hat{a}(y))=(N(y)\bar{b}(\Pi(y)),\bar{a}(\Pi(y)))$ for $y \in [0,2[^d$, with $\Pi(y)=(\pi(y_1),\dots,\pi(y_d))$,
 $\pi(z)=z {\mathbf 1}_{[0,1]}(z) + (2-z) {\mathbf 1}_{]1,2]}(z)$ and $N(y)= {\rm diag}({\rm sign}(1-y_1),\dots,{\rm sign}(1-y_d))$ (${\rm diag}(x)$
denotes the diagonal matrix whose diagonal is given by $x$). Due to the boundedness of $\hat{b}$ and to the continuity and to the ellipticity of $\hat{a}$,
such a martingale problem is uniquely solvable (see Stroock and Varadhan \cite{stroock varadhan 1979}, Chapter VII).
Extend $\Pi$ from a $2{\mathbb Z}^d$-periodicity argument to the whole set $\bR^d$, and derive from
the It\^o-Tanaka formula that $\Pi(\hat{X})$ satisfies the martingale problem associated to $(\bar{b},\bar{a})$ with normal reflection on 
$\partial [0,1]^d$. In particular, $\Pi(\hat{X})$ and $X$ have the same law. It is then sufficient to prove that almost every trajectory of 
$\Pi(\hat{X})$ satisfies the conditions
given in the above paragraph, and thus to prove that almost every trajectory of $\hat{X}$ hits $F \equiv \cup_{k \in 2\bZ^d} (k+F_0)$ and that every point
of the boundary of $F$ is regular for $\hat{X}$.
\vspace{5pt}
\\
Due to the boundedness of $\hat{b}$ and $\hat{a}$ and to the uniform ellipticity of $\hat{a}$, deduce from 
Pardoux \cite{pardoux 1999}, Section 2, that the process $\hat{X}$, seen as a diffusion process with values in the torus $\bbR^d/(2\bbZ^d)$ is recurrent, and
from Pinsky \cite{pinsky 1995}, Theorem 3.3, Section 3, Chapter II, that every point of the boundary of $F$
is regular for $\hat{X}$.
\vspace{5pt}
\\
The second assertion in Corollary \ref{COR1} follows again from the mapping theorem. $\square$
\mysection{Sketch of the Proof of Theorem \ref{MAIN_THM2}}
\label{behaviour}
We now present several ideas to establish Theorem \ref{MAIN_THM2}. Again, the whole proof is given in Delarue \cite{delarue:preprint}.
\vspace{5pt}
\\
We denote by {\bf (\ovA)} the assumption satisfied by $\bar{a}$ and $\bar{b}$ under the statement of Theorem \ref{MAIN_THM2}. Briefly,
$\bar{a}$ and $\bar{b}$ are bounded and Lipschitz continuous, and $\bar{a}$ is uniformly elliptic.
\subsection{Description of the Method}
\label{description}
As a starting point of the proof, recall that the two-dimensional Brownian motion $B$ never hits zero at a positive time, but hits infinitely often
any neighbourhood of zero with probability one. 
The proof of this result (see e.g. Friedman \cite{friedman 1975}) relies on the differential form of the Bessel process of index 1, {\it i.e.} of the process 
$|B|$. In short, for $B$ different from zero, $d|B|$ writes $d|B_t| = 1/(2|B_t|) dt + d \tilde{B}_t$, where $\tilde{B}$ denotes a one-dimensional Brownian motion.
\vspace{5pt}
\\
The common strategy to investigate the recurrence and transience properties of $|B|$ then consists in exhibiting a Lyapunov function for the process 
$|B|$ (see again Friedman \cite{friedman 1975} for a complete
review on this topic). In dimension two, {\it i.e.} in our specific setting, the function \emph{ln} is harmonic for the process $|B|$ (the It\^o formula yields
for $B$ different from zero: $d\ln(|B_t|) = |B_t|^{-1} d \tilde{B}_t$). In particular, for $|B_0|=1$, for a large $N \geq 1$ and for the stopping time 
$\bar{\tau}_N \equiv \tau_{N^{-1}} \wedge \tau_2$, with $\tau_x \equiv \inf\{t \geq 0, |B_t| =x \}$:
\begin{equation*}
\begin{split}
0 = {\mathbb E}(\ln(|B_0|)) &={\mathbb E}(\ln(|B_{\bar{\tau}_N}|))
\\
&= -\ln(N) {\mathbb P}\{\tau_{N^{-1}}<\tau_2 \} + \ln(2) {\mathbb P}\{\tau_2 < \tau_{N^{-1}}\}.
\end{split}
\end{equation*}
Hence, $\lim_{N \rightarrow + \infty} {\mathbb P}\{\tau_{N^{-1}} < \tau_2\} = 0$, so that ${\mathbb P}\{\tau_0<\tau_2\}=0$. Of course, for every $r \geq 2$,
the same property holds: ${\mathbb P}\{\tau_0<\tau_r\}=0$. Letting $r$ tend towards $+ \infty$, it comes ${\mathbb P}\{\tau_0<+\infty\}=0$, as announced above.
\vspace{5pt}
\\
More generally, if $Y$ denotes a Bessel process of index $\gamma \geq -1$, $\gamma \not = 1$, {\it i.e.} $Y$ satisfies the SDE $dY_t = (\gamma/2) Y_t^{-1} + d\tilde{B}_t$, where
$\tilde{B}$ denotes a one-dimensional Brownian motion, the function $\phi(x)=x^{1-\gamma}$ is harmonic for the process $Y$. A similar argument to the
previous one permits to derive the attainability or non-attainability of the origin for the process $Y$, depending on the sign of $1-\gamma$.
\vspace{5pt}
\\
Roughly speaking, the strategy used in Delarue \cite{delarue:preprint} to establish Theorem \ref{MAIN_THM2} aims to reduce the original limit absorption problem 
to the attainability problem for a suitable process. For this reason, 
we translate the absorption problem in $F_0$ to the neighbourhood of zero. This amounts to investigate the
behaviour of ${\mathbb E}(\check{T}_{\ell})$, for $\ell \in ]0,1]$, with $X_0 = (1,1)^t$ and $\check{T}_{\ell} \equiv \inf\{t \geq 0, X_t^1+X_t^2
\leq \ell\}$. As guessed by the reader, our final result then writes in terms of $\bar{a}(0)$ and not of $\bar{a}(1,1)$ as in Theorem \ref{MAIN_THM2}.
\vspace{5pt}
\\
A crucial point is then to define the analogue of $|B|$ for the reflected diffusion $X$. In order to take into account the influence of the diffusion
matrix $\bar{a}(0)$, we focus on the norm of $X$ with respect to the scalar product induced by the inverse of the matrix $\bar{a}(0)$ (which is
non-degenerate under the assumptions of Theorem \ref{MAIN_THM2}):
\begin{equation}
\label{R0}
\begin{split}
\forall t \geq 0, \ Q_t &\equiv \langle X_t, \bar{a}^{-1}(0) X_t \rangle
\\
&= (1-s^2)^{-1} \bigl[ \rho_1^{-2} (X_t^1)^2 + \rho_2^{-2} (X_t^2)^2 - 2 s \rho_1^{-1} \rho_2^{-1} X_t^1 X_t^2 \bigr].
\end{split}
\end{equation}
The process $Q$ is devised to mimic the role played by $|B|^2$ in the non-reflected Brownian case. Its differential form writes:
\begin{prop}
\label{cor9.0}
There exist a constant $C_{\ref{cor9.0}}$, depending only on known parameters in Assumption {\bf (\ovA)}, as well as a measurable function 
$\Gamma_{\ref{cor9.0}} : \bbR^2 \rightarrow \bbR$, bounded by $C_{\ref{cor9.0}}$, such that for $N \geq 1$ and $t \in [0,\zeta_N]$, $\zeta_N \equiv \inf \{ t\geq 0, \, |X_t| \leq N^{-1}\}$:
\begin{equation}
\label{R00}
\begin{split}
\forall t \in [0,\zeta_N], \ dQ_t^{1/2} &=   \Gamma_{\ref{cor9.0}}(X_t) dt + \frac{1}{2} Q_t^{-1/2} dt
\\
&\hspace{15pt} - \frac{s}{\sqrt{1-s^2}}  \bigl[ \rho_1^{-1} dH_t^1 + \rho_2^{-1} dH_t^2 \bigr] -  Q_t^{-1/2} \bigl[ \kappa^1_t dK_t^1 + \kappa^2_t dK_t^2 \bigr]
\\
&\hspace{15pt}+ Q_t^{-1/2} \langle  \bar{\sigma}(X_t) \bar{a}^{-1}(0) X_t,  dB_t \rangle,
\end{split}
\end{equation}
with:
\begin{equation}
\label{Lambda00}
\kappa^1_t \equiv \frac{1}{1-s^2} \bigl[ \rho_1^{-2} - s \rho_1^{-1} \rho_2^{-1} X_t^2  \bigr], \ \kappa^2_t \equiv \frac{1}{1-s^2} \bigl[ \rho_2^{-2} - s \rho_1^{-1} \rho_2^{-1} X_t^1  \bigr].
\end{equation}
\end{prop}
The proof of Proposition \ref{cor9.0} is put aside for the moment.
\vspace{5pt}
\\
At this stage of the sketch, the most simple situation to focus on seems to be the case ``$s=0$''. In this framework, the $dH$ term in \eqref{R00} 
reduces to zero so that $Q$ behaves
like a standard It\^o process in the neighbourhood of zero:
\begin{enumerate}
\item[(1)] For $X$ close to $0$, the scalar product $Q_t^{-1/2} \langle \bar{\sigma}(X_t) \bar{a}^{-1}(0) X_t,dB_t \rangle$ driving the martingale part
of $Q^{1/2}$ looks like the process $Q_t^{-1/2} \langle \bar{\sigma}^{-1}(0)X_t,dB_t \rangle$, which is, thanks to L\'evy's theorem, a Brownian motion.
\item[(2)] For $X$ close to $0$, the term $\Gamma_{\ref{cor9.0}}(X_t)$ is negligible in front of $Q_t^{-1/2}$.
\end{enumerate}
Thus, from (1) and (2), the differential form of the process $Q_t^{1/2}$ satisfies, for $s=0$, the Bessel equation of index 1, at least for $X$ in the
neighbourhood of zero. 
\vspace{5pt}
\\
For the sake of simplicity, we assume from now on that $\overline{a}(x)$ reduces for all $x \in [0,1]^2$ to $\overline{a}(0)$,
and in the same way, that $\Gamma_{\ref{cor9.0}}(x)$ vanishes for $x \in [0,1]^2$. Up to the differential term $dK$, $Q^{1/2}$ then appears as a Bessel
process. For this reason, we expect the function \emph{ln} to be a kind of Lyapunov function (in a sense to be precised) for $Q_t^{1/2}$ in the case $s=0$.
It\^o's formula applied to \eqref{R00} yields for $Q$ different from zero:
\begin{equation}
\label{ITOSKETCH}
d \ln(Q_t^{1/2}) = d \tilde{B}_t - Q_t^{-1} \bigl[ \kappa^1_t dK_t^1 + \kappa^2_t dK_t^2 \bigr].
\end{equation}
It now remains to handle the $dK$ term in \eqref{ITOSKETCH}. Since $s$ vanishes in this first analysis, the processes 
$\kappa^1$ and $\kappa^2$ driven by the differential
terms $dK^1$ and $dK^2$ are always bounded from above and from below by positive constants.
Recall moreover that $Q_t^{-1} dK_t^1 = Q_t^{-1} {\mathbf 1}_{\{X_t^1=1\}} dK_t^1$ (the same holds for $dK^2_t$).
We deduce that there exists a constant $C>0$ such that:
\begin{equation}
\label{SKETCH1}
d \tilde{B}_t - C [dK_t^1 + dK_t^2 ] \leq d \ln(Q_t^{1/2}) \leq d\tilde{B}_t - C^{-1} [ dK_t^1 + dK_t^2 ].
\end{equation}
The strategy then consists in comparing \eqref{SKETCH1} to a bounded real functional of $X$. The most simple one is $|X|^2$. For $\bar{b}=0$, it is plain to derive
from It\^o's formula (note again that $\langle X_t,dH_t \rangle = X_t^1 dH_t^1 + X_t^2 dH_t^2$ vanishes since $dH_t^1 = {\mathbf 1}_{\{X_t^1=0\}} dH_t^1$
and $dH_t^2 = {\mathbf 1}_{\{X_t^2=0\}}dH_t^2$):
\begin{equation}
\label{SKETCH2}
\begin{split}
d|X|_t^2 &= {\rm trace}(\bar{a}(0)) dt - 2 ( X_t^1 dK_t^1 + X_t^2 dK_t^2 ) + 2 \langle X_t,\bar{\sigma}(0)dB_t \rangle
\\
&= {\rm trace}(\bar{a}(0)) dt - 2 (dK_t^1 + dK_t^2) + 2 \langle X_t,\bar{\sigma}(0)dB_t \rangle.
\end{split}
\end{equation}
Add \eqref{SKETCH1} and \eqref{SKETCH2} (up to multiplicative constants) and deduce:
\begin{equation}
\label{SKETCH3}
\begin{split}
&d \bigl( -\ln(Q_t^{1/2}) + \frac{C}{2} |X_t|^2 \bigr) \leq \frac{C}{2} {\rm trace}(\bar{a}(0)) dt - d\tilde{B}_t + C \langle X_t,\bar{\sigma}(0)dB_t \rangle
\\
&d \bigl( - \ln(Q_t^{1/2}) + \frac{C^{-1}}{2} |X_t|^2 \bigr) \geq \frac{C^{-1}}{2} {\rm trace}(\bar{a}(0)) dt - d\tilde{B}_t + C^{-1} \langle X_t,\bar{\sigma}(0)dB_t \rangle.
\end{split}
\end{equation}
Recall that $X_0=(1,1)^t$. Take the expectation between 0 and $\tau_N \equiv \inf \{t \geq 0, \ Q_t^{1/2}=1/N\}$, for $N \geq 1$. Since the expectation of the martingale
parts vanishes and since the matrix $\bar{a}(0)$ is non-degenerate, it comes for a new positive constant $C'$:
\begin{equation*}
(C')^{-1} (\ln(N)-1) \leq {\mathbb E}(\tau_N) \leq C'(\ln(N)+1).
\end{equation*}
Of course, ${\mathbb E}(\tau_N)$ does match exactly ${\mathbb E}(\check{T}_{1/N})$. However, it is well seen from the definition of $Q^{1/2}$ that there exists a constant $c>0$, such that
$\tau_{c^{-1}N} \leq \check{T}_{1/N} \leq \tau_{cN}$. This proves the third point in Theorem \ref{MAIN_THM2}. $\square$
\subsection{Non-Zero Cases}
\label{non_zero_cases}
As easily guessed by the reader, the cases ``$s<0$'' and ``$s>0$'' are more difficult. We derive first several straightforward consequences from 
Proposition \ref{cor9.0}:
\begin{enumerate}
\label{explications2}
\item[(3)] If $s>0$, the $dH$ term is always nonpositive. Up to a slight modification of the $dK$ term (that does not play any role in the neighbourhood of the origin),
the differential form of $Q^{1/2}$ writes as the differential form of $Q^{1/2}$ in
the zero case plus a nonincreasing process. This explains why the process $Q^{1/2}$ reaches the neighbourhood of the origin in the case ``$s>0$'' faster 
than in the case ``$s=0$''. We thus expect the case ``$s>0$'' to be sub-logarithmic (again, in a sense to be precised).
\item[(4)] On the opposite, if $s<0$, the $dH$ term is always nonnegative. With a similar argument to the previous one, we expect the case ``$s<0$'' to be
super-logarithmic.
\end{enumerate}
The strategy now consists in correcting $Q^{1/2}$ to get rid of the $dH$ term and to reduce the analysis to the one stated for ``$s=0$''.
It is then rather natural to focus on the auxiliary process $\forall t \geq 0, \ A_t \equiv Q_t^{1/2} + Z_t$, with:
\begin{equation*}
\forall t \geq 0, \ Z_t \equiv \frac{s}{\sqrt{1-s^2}} \langle (\rho_1^{-1},\rho_2^{-1})^t,X_t \rangle = \frac{s}{\sqrt{1-s^2}} \bigl[ \rho_1^{-1} X_t^1 + \rho_2^{-1} X_t^2 \bigr].
\end{equation*}
For $s$ close to 0, the processes $A$ and $Q^{1/2}$ are equivalent in the following sense: there exists a constant $C>0$, such that 
$C^{-1} Q_t^{1/2} \leq A_t \leq C Q_t^{1/2}$. The value of the constant $C$ can be precised for $s$ in the neighbourhood of zero.
Indeed, the process $Z$ can be expressed as: $Z_t = s \gamma_t(s) Q_t^{1/2}$, with $1-\varepsilon \leq  \gamma_t(s)
 \leq \sqrt{2}+ \varepsilon$, for $s$ in the neighbourhood of zero.
\vspace{5pt}
\\
Actually, we prove in Delarue \cite{delarue:preprint} that the equivalence property still holds true for all $s \in ]-1,1[$. As a first consequence, 
we derive that the processes $Q^{1/2}$ and $A$ admit the same asymptotic behaviour. We thus focus on the second one and in particular on the differential
form of $A$ (the proof derives from Proposition \ref{cor9.0} and is thus left to the reader):
\begin{prop}
\label{prop9.4}
There exist a constant $C_{\ref{prop9.4}}$, depending only on known parameters in {\bf (\ovA)}, and a function $\Gamma_{\ref{prop9.4}}$, bounded by $C_{\ref{prop9.4}}$,
such that:
\begin{equation}
\label{R2A}
\begin{split}
\forall t \in [0,\zeta_N], \ dA_t &=   \Gamma_{\ref{prop9.4}}(X_t) dt + \frac{1}{2} Q_t^{-1/2} dt
\\
&\hspace{15pt} - \bigl[ \bar{\kappa}^1_t dK_t^1 + \bar{\kappa}^2_t dK_t^2 \bigr]
\\
&\hspace{15pt}+ \langle  Q_t^{-1/2} \bar{\sigma}(X_t) \bar{a}^{-1}(0) X_t + \frac{s}{\sqrt{1-s^2}} \bar{\sigma}(X_t) (\rho_1^{-1},\rho_2^{-1})^t,  dB_t \rangle,
\end{split}
\end{equation}
with:
\begin{equation*}
\begin{split} 
\forall t \geq 0, \ \bar{\kappa}^1_t &=  Q_t^{-1/2} \kappa^1_t + \frac{s}{\sqrt{1-s^2}} \rho_1^{-1},
\\
\bar{\kappa}^2_t &=  Q_t^{-1/2} \kappa^2_t + \frac{s}{\sqrt{1-s^2}} \rho_2^{-1}.
\end{split}
\end{equation*}
\end{prop}
From now on, we assume that $\bar{\sigma}$ and $\bar{a}$ do not depend on $x$, as done in the previous subsection 
to investigate the case ``$s=0$''.
In the same way, we assume for the sake of simplicity that $\Gamma_{\ref{prop9.4}}$ reduces to zero. Equation
\eqref{R2A} becomes:
\begin{equation}
\label{R2B}
\begin{split}
\forall t \in [0,\zeta_N], \ dA_t &=   \frac{1+s \gamma_t(s)}{2} A_t^{-1} dt
\\
&\hspace{15pt} - \bigl[ \bar{\kappa}^1_t dK_t^1 + \bar{\kappa}^2_t dK_t^2 \bigr]
\\
&\hspace{15pt}+ \langle  Q_t^{-1/2} \bar{\sigma}^{-1}(0) X_t + \frac{s}{\sqrt{1-s^2}} \bar{\sigma}(0) (\rho_1^{-1},\rho_2^{-1})^t,  dB_t \rangle,
\end{split}
\end{equation}
We try to mimic the arguments given in Subsection \ref{description}. For $s$ close to zero, we deduce from Proposition \ref{cor9.0}
that processes $\kappa^1$ and $\kappa^2$ are positive. Hence, from Proposition \ref{prop9.4}, the same holds for $\bar{\kappa}^1$ and $\bar{\kappa}^2$.
Since the general case $s \in ]-1,1[$ is really more difficult to handle, we restrict our analysis to $s$ in the neighbourhood of zero.
\vspace{5pt}
\\
It now remains to investigate the martingale part in \eqref{R2A}. Its bracket writes:
\begin{equation*}
\begin{split}
\frac{d \langle A \rangle_t}{dt} &= \bigl| Q_t^{-1/2} \bar{\sigma}^{-1}(0) X_t + \frac{s}{\sqrt{1-s^2}} \bar{\sigma}(0)(\rho_1^{-1},\rho_2^{-1})^t \bigr|^2
\\
&=\bigl[ 1 + 2 Q_t^{-1/2} Z_t
+ \frac{s^2}{1-s^2} \langle (\rho_1^{-1},\rho_2^{-1})^t,\bar{a}(0) (\rho_1^{-1},\rho_2^{-1})^t \rangle  \bigr]
\\
&=\bigl[ 1 + 2s \gamma_t(s) + 2s^2 (1-s)^{-1}  \bigr]
\end{split}
\end{equation*}
Apply now the function $x^{\beta}$ to $A$, for a given $\beta \not =0$ (recall that this function is harmonic for the Bessel process of index $1-\beta$).
It\^o's formula yields:
\begin{equation*}
\begin{split}
d A_t^{\beta} &= \frac{\beta}{2} A_t^{\beta-2} \Theta_t(\beta,s)dt - \beta A_t^{\beta-1} \bigl[ \bar{\kappa}^1_t dK_t^1 + \bar{\kappa}^2_t dK_t^2 \bigr] + dM_t,
\\
&{\rm with} \ \Theta_t(\beta,s) \equiv \bigl[ (\beta-1)(1 + 2s \gamma_t(s) + 2s^2(1-s)^{-1}) + 1 + s \gamma_t(s) \bigr]
\end{split}
\end{equation*}
where $M$ denotes a martingale. 
\vspace{5pt}
\\
We are now in position to complete the proof:
\begin{enumerate}
\item If $s>0$ (and small), then the quantity $\Theta_t(\beta,s)$
tends to $-s \gamma_t(s) - 2s^2(1-s)^{-1}$ when $\beta$ tends to zero. In particular, there exists 
$\beta>0$ ($\beta$ close to zero) such that $\Theta_t(\beta,s) \leq 0$. For this $\beta$, 
$dA_t^{\beta} \leq - \beta A_t^{\beta-1} [ \bar{\kappa}^1_t dK_t^1 + \bar{\kappa}^2_t dK_t^2 ]
+dM_t$. We thus recover the right hand side in \eqref{SKETCH1} with respect to the function $x^{\beta}$ and to the process $A$. We then derive
the second inequality in \eqref{SKETCH3}, with $-\ln(x)$ replaced by $-x^{\beta}$, and $Q^{1/2}$ by $A$.
\vspace{5pt}
\\
This proves the first point in Theorem \ref{MAIN_THM2}.
\item If $s<0$ (and small), we can follow the previous argument and thus prove that there exists $\beta<0$ ($\beta$
close to zero) such that $\Theta_t(\beta,s)\geq 0$. For this $\beta$, $dA_t^{\beta} \leq - \beta A_t^{\beta-1} [ \bar{\kappa}^1_t dK_t^1 + \bar{\kappa}^2_t dK_t^2 ]
+dM_t$. Again, we recover the right hand side in \eqref{SKETCH1} with respect to the function $x^{\beta}$ and to the process $A$, but with
$-C^{-1}>0$, {\it i.e.} $C<0$. We then
derive the second inequality in \eqref{SKETCH3}, with $-\ln(x)$ replaced by $-x^{\beta}$, $Q^{1/2}$ by $A$ and $C$ by $-C$.
\vspace{5pt}
\\
This proves the lower bound in the second point of Theorem \ref{MAIN_THM2} (with respect to a possibly different $\beta_{-}$).
\vspace{5pt}
\\
To obtain the upper bound, note for $\beta \rightarrow - \infty$ that $\Theta_t(\beta,s) \sim \beta (1+ 2 s \gamma_t(s) + 2s^2(1-s)^{-1})$
and deduce (at least for $s$ close to zero) that there exists 
$\beta<0$ such that $\Theta_t(\beta,s)\leq 0$. For this $\beta$, $dA_t^{\beta} \geq - \beta A_t^{\beta-1} [ \bar{\kappa}^1_t dK_t^1 + \bar{\kappa}^2_t dK_t^2 ]
+dM_t$. Again, we follow the argument given in the case ``$s=0$'' and then derive the upper bound in the second point of Theorem \ref{MAIN_THM2} 
(the resulting $\beta_+$ may differ from the one given in the statement).
$\square$
\end{enumerate}
\subsection{Proof of Proposition \ref{cor9.0}}
Write first the SDE satisfied by the process $Q$:
\begin{lem} 
\label{9.2}
Keep the notations introduced in the statement of Proposition \ref{cor9.0}. Then, there exists a constant $C_{\ref{9.2}}$, 
depending only on known parameters in {\bf (\ovA)}, as well as a measurable function $\Gamma_{\ref{9.2}}$, bounded by $C_{\ref{9.2}}$, 
such that:
\begin{equation}
\label{R2}
\begin{split}
\forall t \geq 0, \ dQ_t &= Q_t^{1/2} \Gamma_{\ref{9.2}}(X_t) dt + 2 dt
\\
&\hspace{15pt} - \frac{2s}{\sqrt{1-s^2}}  Q_t^{1/2} \bigl[ \rho_1^{-1} dH_t^1 + \rho_2^{-1} dH_t^2 \bigr] 
- 2 \bigl[ \kappa^1_t dK_t^1 + \kappa^2_t dK_t^2 \bigr]
\\
&\hspace{15pt}+ 2 \langle  \bar{\sigma}(X_t) \bar{a}^{-1}(0) X_t,  dB_t \rangle.
\end{split}
\end{equation}
\end{lem}
Assume for the moment that Lemma \ref{9.2} holds and complete the proof of Proposition \ref{cor9.0}.
Apply It\^o's formula to $Q^{1/2}$ and derive from the differential form given in \eqref{R2} 
(note that $Q$ does not vanish for $t \in [0,\zeta_N]$):
\begin{equation}
\label{9.3-1}
\begin{split}
\forall t \in [0,\zeta_N], \ dQ_t^{1/2} &=  \frac{1}{2}  \Gamma_{\ref{9.2}}(X_t) dt + Q_t^{-1/2} dt
\\
&\hspace{15pt} - \frac{s}{\sqrt{1-s^2}}  \bigl[ \rho_1^{-1} dH_t^1 + \rho_2^{-1} dH_t^2 \bigr] - Q_t^{-1/2} \bigl[ \kappa^1_t dK_t^1 + \kappa^2_t dK_t^2 \bigr]
\\
&\hspace{15pt} + Q_t^{-1/2} \langle  \bar{\sigma}(X_t) \bar{a}^{-1}(0) X_t,  dB_t \rangle
\\
&\hspace{15pt}  - \frac{1}{2} Q_t^{-3/2} \langle X_t, \bar{a}^{-1}(0) \bar{a}(X_t) \bar{a}^{-1}(0) X_t \rangle dt
\\
&\equiv \bigl[ \Delta(1,X_t)  + Q_t^{-1/2} \bigr] dt + d \Delta_t (2) + d \Delta_t(3) - \frac{1}{2} Q_t^{-3/2} \Delta(4,X_t) dt.
\end{split}
\end{equation} 
Focus on $\Delta(4,\cdot)$ in \eqref{9.3-1}:
\begin{equation}
\label{9.3-2}
\begin{split}
\forall t \in [0,\zeta_N], \ \Delta(4,X_t) &= Q_t  + \langle X_t, \bar{a}^{-1}(0) [ \bar{a}(X_t) - \bar{a}(0) ]
\bar{a}^{-1}(0) X_t \rangle
 \\
  &\equiv Q_t + \Delta(5,X_t).
 \end{split}
 \end{equation}
 Thanks to Assumption {\bf ($\bar{\mathbf A}$)} (boundedness, ellipticity and Lipschitz continuity of $\bar{a}$), 
 the function $\Delta(5,\cdot)$ satisfies:
 \begin{equation}
 \label{9.3-3}
\forall t \in [0,\zeta_N], \ |\Delta(5,X_t)| \leq C Q_t^{3/2}.
 \end{equation}  
where the constant $C$ depends only on known parameters. 
Plug now \eqref{9.3-2} into \eqref{9.3-1}:
\begin{equation}
\label{9.3-4}
\forall t \in [0,\zeta_N], \ dQ_t^{1/2} =  \bigl[\Delta(1,X_t) - \frac{1}{2}  Q_t^{-3/2} \Delta(5,X_t)+ \frac{1}{2} Q_t^{-1/2} \bigr] dt 
 + d \Delta_t (2) + d \Delta_t(3).
\end{equation} 
 Set finally $\displaystyle \Gamma_{\ref{cor9.0}}(X_t) =  \Delta(1,X_t) - \frac{1}{2} Q_t^{-3/2} \Delta(5,X_t)$. 
 Thanks to \eqref{9.3-1}, \eqref{9.3-3} and to \eqref{9.3-4}, this completes
 the proof. $\square$
 \vspace{5pt}
 \\
{\bf Proof of Lemma \ref{9.2}.} Thanks to It\^o's formula, $d\langle X ,\bar{a}^{-1}(0) X \rangle$ writes:
\begin{equation}
\label{R}
\begin{split}
\forall t \geq 0, \ d\langle X_t,\bar{a}^{-1}(0) X_t \rangle &= 2 \langle X_t, \bar{a}^{-1}(0) \bar{b}(X_t) \rangle dt 
+ {\rm trace}\bigl[ \bar{\sigma}^{-1}(0)\bar{a}(X_t) \bar{\sigma}^{-1}(0)\bigr] dt
\\
&\hspace{15pt} + 2 \langle X_t, \bar{a}^{-1}(0) (dH_t - dK_t) \rangle + 
2 \langle  \bar{a}^{-1}(0)X_t, \bar{\sigma}(X_t) dB_t \rangle
\\
&= 2 dt 
\\
&\hspace{15pt} + 2 \langle X_t, \bar{a}^{-1}(0) \bar{b}(X_t) \rangle dt + {\rm trace}\bigl[ \bar{\sigma}^{-1}(0) \bigl( \bar{a}(X_t) - \bar{a}(0) \bigr) \bar{\sigma}^{-1}(0)\bigr] dt
\\
&\hspace{15pt} + 2 \langle X_t, \bar{a}^{-1}(0) (dH_t - dK_t) \rangle + 
2 \langle  \bar{a}^{-1}(0)X_t, \bar{\sigma}(X_t) dB_t \rangle
\end{split}
\end{equation}
Set $\Gamma_{\ref{9.2}}(X_t) = Q_t^{-1/2} \bigl[ 2 \langle X_t, \bar{a}^{-1}(0) \bar{b}(X_t) \rangle  + {\rm trace}[ \bar{\sigma}^{-1}(0) \bigl( \bar{a}(X_t) - \bar{a}(0) \bigr)
\bar{\sigma}^{-1}(0)]\bigr]$.
\vspace{5pt}
\\
Focus now on the $dH$ and $dK$ terms in the above r.h.s. Due to \eqref{LEAST}, note first that $X_t^1 dH_t^1= 0$ and $X_t^2 dH_t^2=0$ and
that $\rho_2^{-1} X_t^2 dH_t^1 = (1-s^2)^{1/2} Q_t^{1/2} dH_t^1$ and  $\rho_1^{-1} X_t^1 dH_t^2 = (1-s^2)^{1/2} Q_t^{1/2} dH_t^2$. Derive that:
\begin{equation}
\label{R,H}
\begin{split}
\langle X_t, \bar{a}^{-1}(0) dH_t \rangle &= X_t^1 \bigl( \bar{a}^{-1}(0) \bigr)_{1,1} dH_t^1 + X_t^2 \bigl( \bar{a}^{-1}(0) \bigr)_{2,2} dH_t^2
\\
&\hspace{15pt} + \bigl( \bar{a}^{-1}(0) \bigr)_{1,2} \bigl[  X_t^2 dH_t^1 + X_t^1 dH_t^2 \bigr]
\\
&\hspace{0pt} = - \frac{s}{1-s^2} \rho_1^{-1} \rho^{-1}_2 \bigl[ X_t^2 dH_t^1 +  X_t^1 dH_t^2 \bigr]
\\
&\hspace{0pt} = - \frac{s}{\sqrt{1-s^2}} Q_t^{1/2} \bigl[\rho_1^{-1} dH_t^1+  \rho_2^{-1}  dH_t^2  \bigr].
\end{split}
\end{equation}
Note again from \eqref{LEAST} that $X_t^1 dK_t^1=dK_t^1$ and $X_t^2dK_t^2=dK_t^2$. Thus:
\begin{equation}
\label{R,K}
\begin{split}
\langle X_t, \bar{a}^{-1}(0) dK_t \rangle &= X_t^1 \bigl( \bar{a}^{-1}(0) \bigr)_{1,1} dK_t^1 + X_t^2 \bigl( \bar{a}^{-1}(0) \bigr)_{2,2} dK_t^2
\\
&\hspace{15pt} + \bigl( \bar{a}^{-1}(0) \bigr)_{1,2} \bigl[ X_t^1 dK_t^2 + X_t^2 dK_t^1 \bigr]
\\
&\hspace{0pt} = \frac{1}{1-s^2} \bigl[ \rho_1^{-2} - s \rho_1^{-1} \rho_2^{-1} X_t^2 \bigr] dK_t^1+ \frac{1}{1-s^2} \bigl[
\rho_2^{-2} - s \rho_1^{-1} \rho_2^{-1} X_t^1 \bigr] d K_t^2.
\end{split}
\end{equation}
Plug now \eqref{R,H} and \eqref{R,K} into \eqref{R} and complete the proof of
\eqref{R2}. $\square$
\vspace{5pt}
\\
Note from \eqref{9.3-1}, \eqref{9.3-2}, \eqref{9.3-4} and \eqref{R} that
$\Gamma_{\ref{cor9.0}}$ reduces to zero if $\bar{b}$ vanishes and $\bar{a}$ is constant.
This explains why we can assume in Subsections \ref{description} and \ref{non_zero_cases} $\bar{a}(x)=\bar{a}(0)$, $\bar{b}=0$ and $
\Gamma_{\ref{cor9.0}} = \Gamma_{\ref{prop9.4}}=0$.
\mysection{Conclusion}
\label{Conclusion}
Probabilistic models depending on an extra random environment have been widely developed for twenty years. From a practical
point of view, the environment expresses external constraints that satisfy some statistical properties. For example, the behaviour 
of an investor that borrows a certain amount of money may depend on the values of various assets as well as some quasi-periodic
parameters.
\vspace{5pt}
\\
From a more theoretical point of view, our work shows that earlier results (see e.g. Louchard \cite{louchard 1995}, Louchard et al.
\cite{louchard schott 1991} and \cite{louchard tolley zimmerman 1994} and Guillotin and Schott \cite{guillotin schott 2002}) for the diffusive behaviour of the banker algorithm 
extend to a more general framework. This is a new small step in the robustness of the underlying model.
The theory of stochastic homogenization is also relevant for other distributed algorithms such as the problem of
colliding stacks.
\vspace{5pt}
\\
The reader may object that our own paper exclusively focuses on the trend-free case (see Assumption {\bf (A.4)}), whereas so-called contracting and expanding regimes for
distributed algorithms have been investigated in the previous references. We reasonably guess that the strategy presented here applies
to the expanding case, after a suitable substraction of the trend, as done for example in Louchard and Schott \cite{louchard schott 1991}. On the opposite, the
contracting framework would require further developments, both based on homogenization techniques and on large deviations arguments as in 
Maier \cite{maier 1991}. We leave this for further research.
\vspace{5pt}
\\
In a different perspective, two questions follow from the study of the limit reflected diffusion. On the one hand, our global result writes as
a double asymptotic property: we let first $m$ tend to $+ \infty$ in $m^{-2} T^{(m)}$ and then $\lambda$ to 2 in ${\mathbb E}(T_\lambda)$ (at least in
the two dimensional setting). A natural question would consist in investigating $m^{-2} {\mathbb E}(T^{(m)}_{\Lambda_m})$ (with obvious notations)
as $m$ tends to $ + \infty$ and for $(m^{-1} \Lambda_m)_{m \geq 1}$ converging towards 2. On the other hand, the extension of Theorem  \ref{MAIN_THM2} to the $d$-dimensional
case appears as a challenging problem. Nevertheless, numerical methods offer a possible alternative to obtain an empirical description of the asymptotic
behaviour of the limit diffusion or even to compute with a Monte-Carlo procedure an estimate of underlying quantities of interest.


\begin{thebibliography}{99}
\bibitem {blp 1978} 
Bensoussan, A., Lions, J.L., Papanicolaou, G. (1978): Asymptotic analysis for periodic structures. North-Holland, Amsterdam.

\bibitem{billingsley 1968}
Billingsley, P. (1968): Convergence of probability measures. John Wiley \& Sons, Inc., New York-London-Sydney.


\bibitem{dacunha-castelle duflo t2} 
Dacunha-Castelle, D., Duflo, M. (1983): Probabilit\'es et statistiques, Tome 2, Probl\`emes \`a temps mobile. First Edition.
Masson, Paris.

\bibitem{delarue:preprint}
Delarue, F. (2005): Recurrence and Transience of a Reflected Diffusion in the Square. Preprint of the LPMA, Universities Paris VI and VII.

\bibitem{ellis 1977}
Ellis, C. A. (1977):
Probabilistic models of computer deadlock.
Inform. Sci. 12, 43--60.


\bibitem{flajolet 1986}
Flajolet, P. (1986): The Evolution of Two Stacks in Bounded Space and Random Walks in a Triangle.
Proceedings of FCT'86., LNCS, No 233, 425--340, Springer Verlag.

\bibitem{freidlin 1964}
Freidlin, M., (1964): The Dirichlet problem for an equation with periodic coefficients depending on a small parameter. Teor. Veroyatnost. I. Primenen 9, 133-139.

\bibitem{friedman 1975}
Friedman, A. (1975): Stochastic differential equations and applications. Vol. 1. 
Probability and Mathematical Statistics, Vol. 28. Academic Press New York-London.

\bibitem{gilbarg and trudinger 1983}
Gilbarg, D., Trudinger, N.S. (1983): 
Elliptic partial differential equations of second order, Second edition, Grundlehren der Mathematischen Wissenschaften, 224, Springer-Verlag, Berlin, 1983.

\bibitem{guillotin schott 2002}
Guillotin-Plantard, N., Schott, R. (2002): Distributed algorithms with dynamic random transitions. Random Structures and Algorithms, 21, 3-4, 371-396. 

\bibitem{haberman 1978}
Haberman ~A.~N. (1978):  System Deadlocks. K.~M. Chandy and R.~T. Yeh, Vol. 3, 256--297, Prentice-Hall. 

\bibitem{jacod shiryaev 1987}
Jacod, J., Shiryaev A.N. (2003): Limit theorems for stochatic processes. Grundlehren der mathematischen Wissenschaften, Springer-Verlag, Berlin, Heidelberg.

\bibitem{jko 1994}
Jikov, V.V., Kozlov, S.M., Oleinik, O.A. (1994): Homogenization of differential operators and integral functionals. Springer-Verlag, Berlin.

\bibitem{krylov 1979}
Krylov, N.V (1979): Controlled diffusion processes. Springer-Verlag, New-York.

\bibitem{knuth 1973}
Knuth~D.~E. (1973):  The Art of Computer Programming, Vol.1, Addison-Wesley.  

\bibitem{lions sznitman 1984}
Lions, P.-L., Sznitman, A.-S. (1984):
Stochastic differential equations with reflecting boundary conditions.
Comm. Pure Appl. Math., 37, 511--537.

\bibitem{louchard 1995}
Louchard, G. (1995): Some Distributed Algorithms Revisited. 
Commun. Statist.- Stochastic models, Vol. 11, No 4, 563--586.

\bibitem{louchard schott 1991}
Louchard, G., Schott, R. (1991): Probabilistic Analysis of Some Distributed Algorithms. 
Random Structures and Algorithms, No 2, 151--186.

\bibitem{louchard tolley zimmerman 1994}
Louchard~G., Schott~R., Tolley~M. and Zimmermann~P. (1994):
Random Walks, Heat Equations and Distributed Algorithms. Computational and Applied Mathematics, 53, 243-274.

\bibitem{maier 1991}
Maier~R.~S. (1991):  Colliding Stacks: A Large Deviations Analysis. 
Random Structures and Algorithms , No 2, 379--420.

\bibitem{maier schott 1993}
Maier~R.~S., Schott~R. (1993): Exhaustion of Shared Memory: Stochastic Results. Proceedings of WADS'93,
LNCS No 709, 494--505, Springer Verlag.

\bibitem{olla 2003}
Olla, S. (2001):  
Central limit theorems for tagged particles and for diffusions in random environment. 
Comets, Francis (ed.) et al., Random media. Paris: Soci\'et\'e Math\'ematique de France. Panor. Synth. 12, 75-100.


\bibitem{pardoux 1999}
Pardoux, \'E (1999): Homogenization of linear and semilinear second order parabolic PDEs with periodic coefficients: a probabilistic approach. 
J. Funct. Anal.,  167 ,498--520.

\bibitem{pardoux veretennikov 1997}
Pardoux, \'E., Veretennikov, A. Yu. (1997): 
Averaging of backward stochastic differential equations, with application to semi-linear PDE's. Stochastics Stochastics Rep., 60, 255-270.

\bibitem{pinsky 1995}
Pinsky, R.G (1995): Positive harmonic functions and diffusion. Cambridge Studies in Advanced Mathematics, 45, Cambridge University Press, Cambridge, 1995.

\bibitem{rogers williams T1}
Rogers, L.C.G., Williams, D. (1979): Diffusions, Markov Proceses and Martingales: Volume 1, Foundations. John Wiley and Sons Ltd, Chichester, UK.

\bibitem{rogers williams T2}
Rogers, L.C.G., Williams, D. (1987): Diffusions, Markov Proceses and Martingales: Volume 2, It\^o Calculs. John Wiley and Sons Ltd, Chichester, UK.

\bibitem{saisho 1987}
Saisho, Y. (1987): Stochastic differential equations for multidimensional domain with reflecting boundary.  
Probab. Theory Related Fields,  74, 455--477. 


\bibitem{slominski 2001} 
Slomi\'nski, L. (2001):
Euler's approximations of solutions of SDEs with reflecting boundary. 
Stochastic Process. Appl., 94, 317--337.

\bibitem{stroock varadhan 1979} Stroock, D.W., Varadhan S.R.S, 1(979): Multidimensional diffusion processes.
Springer-Verlag, New-York.

\bibitem{tanaka 1979}
Tanaka, H. (1979):
Stochastic differential equations with reflecting boundary condition in convex regions.  Hiroshima Math. J.,  9, 
163--177. 

\bibitem{yao 1981}
Yao~A.~C. (1981): An Analysis of a Memory Allocation Scheme for Implementing Stacks. 
SIAM J.~Comput., No 10, 398--403.

\end{thebibliography}
\end{document}